\documentclass[11pt]{article}
\usepackage[sectionbib]{natbib}
\usepackage{array,epsfig,fancyheadings,rotating}
\usepackage[colorlinks,linkcolor=blue,anchorcolor=blue,citecolor=blue,CJKbookmarks=True]{hyperref}
\usepackage{sectsty, secdot}
\usepackage{amsbsy}
\usepackage{amsmath}
\usepackage{amssymb}
\usepackage{epsfig}
\usepackage{multicol}
\usepackage{multirow}
\usepackage{color}
\usepackage{bm}
\usepackage{threeparttable}
\usepackage{float}
\usepackage{array}
\usepackage{amsmath}
\usepackage{amssymb}
\usepackage{amsfonts}
\usepackage{multirow}
\usepackage{amsthm}
\usepackage{graphicx}
\usepackage{bm}
\usepackage{enumerate}
\usepackage{subfigure}
\usepackage{booktabs}
\usepackage{pdflscape}
\usepackage{url}
\usepackage{xcolor,graphicx,float}
\usepackage{rotating}
\usepackage{times}
\usepackage{indentfirst}
\usepackage{blindtext}
\usepackage{algorithm,algorithmicx,float}
\usepackage{lipsum}
\usepackage[noend]{algpseudocode}
\usepackage{setspace}
\usepackage{booktabs}

\usepackage[sectionbib]{natbib}
\usepackage{array,epsfig,fancyheadings,rotating}
\usepackage[]{hyperref}  %<----modified by Ivan
\usepackage{sectsty, secdot}
\usepackage{authblk}
\usepackage{amsmath}
\usepackage{amssymb}
\usepackage{amsfonts}
\usepackage{multirow}
\usepackage{amsthm}
\usepackage{graphicx}
\usepackage{bm}
\usepackage{enumerate}
\usepackage{subfigure}
\usepackage{booktabs}
\usepackage{pdflscape}
\usepackage{url}
\usepackage{xcolor,graphicx,float}
\usepackage{rotating}
\usepackage{times}
\usepackage{indentfirst}
\usepackage{blindtext}
\usepackage{algorithm,algorithmicx,float}
\usepackage{lipsum}
\usepackage[noend]{algpseudocode}
\usepackage{setspace}
\usepackage{booktabs}


\makeatother
\renewcommand{\baselinestretch} {1.5}
\makeatletter 
\def\singlespace{\def\baselinestretch{1}\@normalsize}

\newlength\savewidth

\@addtoreset{equation}{section}
\newtheorem{theorem}{Theorem}

\newtheorem{definition}{Definition}
\newtheorem{assumption}{Assumption}
\newtheorem{remark}{Remark}%[section]
\newtheorem{corollary}{Corollary}%[section]
\newtheorem{example}{Example}%[section]
\newcommand{\by}{\mathbf{y}}
\newcommand{\bx}{\mathbf{x}}
\newcommand{\bk}{\mathbf{k}}
\newcommand{\bw}{\mathbf{w}}
\newcommand{\bX}{\mathbf{X}}
\newcommand{\bP}{\mathbf{P}}
\newcommand{\bD}{\mathbf{D}}
\newcommand{\bI}{\mathbf{I}}
\newcommand{\bmu}{\boldsymbol{\mu}}
\newcommand{\bbeta}{\boldsymbol{\beta}}
\newcommand{\bgamma}{\boldsymbol{\gamma}}
\newcommand{\blambda}{\boldsymbol{\lambda}}
\renewcommand{\hat}{\widehat}
\def\tilde{\widetilde}
\DeclareMathOperator*{\argmin}{arg\,min}
\DeclareMathOperator*{\argmax}{arg\,max}

\UseRawInputEncoding
\def\tilde{\widetilde}

\def\var{\mbox{var}}
\def\cov{\mbox{cov}}

\newcommand{\bg}{\begin{eqnarray}}
\newcommand{\ed}{\end{eqnarray}}
\newcommand{\bgn}{\begin{eqnarray*}}
\newcommand{\edn}{\end{eqnarray*}}

\makeatletter
\def\singlespace{\def\baselinestretch{1}\@normalsize}

\title{Model averaging: A shrinkage perspective}

\author{Jingfu Peng\thanks{E-mail: aypjf@tsinghua.edu.cn}}
\affil{Yau Mathematical Sciences Center, Tsinghua University}

\date{}

\begin{document}
\begin{sloppypar}

\maketitle

\begin{abstract}

Model averaging (MA), a technique for combining estimators from a set of candidate models, has attracted increasing attention in machine learning and statistics. In the existing literature, there is an implicit understanding that MA can be viewed as a form of shrinkage estimation that draws the response vector towards the subspaces spanned by the candidate models. This paper explores this perspective by establishing connections between MA and shrinkage in a linear regression setting with multiple nested models. We first demonstrate that the optimal MA estimator is the best linear estimator with monotonically non-increasing weights in a Gaussian sequence model. The Mallows MA (MMA), which estimates weights by minimizing the Mallows' $C_p$ over the unit simplex, can be viewed as a variation of the sum of a set of positive-part Stein estimators. Indeed, the latter estimator differs from the MMA only in that its optimization of Mallows' $C_p$ is within a suitably relaxed weight set. Motivated by these connections, we develop a novel MA procedure based on a blockwise Stein estimation. The resulting Stein-type MA estimator is asymptotically optimal across a broad parameter space when the variance is known. Numerical results support our theoretical findings. The connections established in this paper may open up new avenues for investigating MA from different perspectives. A discussion on some topics for future research concludes the paper.

\end{abstract}


\section{Introduction}

Model averaging (MA) is an umbrella term for methods that combine multiple candidate models to make a decision, typically in regression and forecasting problems. The concept of MA was first introduced by \cite{Laplace1818deuxieme} \citep[see][for a comprehensive review of the historical development of MA]{Moral-Benito2015}. In recent years, MA has received an explosion of interest in both machine learning and statistics \citep[see, e.g.,][]{Fletcher2018}. It is regarded as a viable alternative to model selection (MS) techniques, as it aims to mitigate MS variability and control modeling biases among candidate models. The benefits of MA compared to MS have been theoretically studied in \cite{Peng2022improvability}.

In the existing literature, MA has been approached using either Bayesian or frequentist frameworks. The Bayesian perspective on MA can be found in works such as \cite{Draper1995}, \cite{George1997}, and \cite{Hoeting1999}. Within the frequentist paradigm, MA strategies have become increasingly popular in forecasting literature since the works of \cite{Barnard1963} and \cite{Bates1969} \citep[see][for a review on the literature]{Timmermann2006Forecast}. In recent years, several important techniques has been developed, including boosting \citep{FREUND1995256}, bagging \citep{Breiman1996b}, stacking \citep{Wolpert1992Stacked, Breiman1996stacking}, random forest \citep{Amit1997}, information criterion weighting \citep{Buckland1997, Hjort2003}, adaptive regression by mixing \citep{Yang2001Adaptive, Yuan2005}, exponential weighting \citep{George1986, Leung2006, Dalalyan2012Sharp}, and Q-aggregation \citep{Rigollet2012kl, Dai2014, Bellec2018}. Additionally, there is a growing body of literature focused on constructing asymptotically optimal MA, with the goal of finding the optimal convex combination of candidate estimates. This is typically achieved by minimizing performance measures such as Mallows' $C_p$ \citep[see, e.g.,][]{Blaker1999adaptive, Hansen2007least, Wan2010least} or cross-validation error \citep[see, e.g.,][]{Hansen2012Jackknife, Zhang2013jackknife, Ando2014high}.

Despite the extensive theoretical work and wide applications of MA, there is a commonly held viewpoint that MA is essentially a shrinkage estimator, and that other shrinkage methods can also achieve the objectives of MA. This view has been substantiated by several studies. For instance, the results in Section 5.1 of \cite{kneip1994ordered} indicate that combining two linear smoothers by minimizing Mallows' $C_p$ yields a James-Stein estimator. The relationship between Mallows model averaging (MMA) and Stein shrinkage estimation has been further explored in \cite{Blaker1999adaptive}, \cite{Hansen2007least}, and \cite{Hansen2014risk} in the context of two nested models. In a semiparametric regression setting, \cite{ullah2017semiparametric} established the connection between MA and ridge shrinkage on the basis of the orthogonal model. Additionally, in a Gaussian location model, \cite{green1991james} proposed a James-Stein type estimator to estimate the best linear combination of two independent biased and unbiased estimators. The methodology in \cite{green1991james} has been further explored by \cite{kim2001james}, \cite{judge2004semiparametric}, and \cite{mittelhammer2005combining}. More recently, \cite{Hansen2016efficient} proposed a Stein method to combine the restricted and unrestricted maximum likelihood estimators in a local asymptotic framework, and showed the asymptotic risk of this shrinkage estimator is strictly less than that of the maximum likelihood estimator.

Note that most of the aforementioned studies focused on the relationship between MA and shrinkage in a two-model setting. The fundamental question that remains to be explored is whether these links persist in the context of multiple models. If so, can some state-of-the-art shrinkage techniques be employed to create MA estimators that perform as optimally as the infeasible optimal convex combinations of candidate models (i.e., the asymptotically optimal MA estimators)? Such answers would have a significant impact on the theories and applications of MA.

This paper addresses the previously mentioned questions in a general linear model setting with multiple nested candidate models. The main contribution is twofold. First, we demonstrate that the optimal MA estimator is equivalent to the optimal linear estimator with monotonically non-increasing weights in a specific Gaussian sequence model. And the MMA estimator \citep{Hansen2007least}, which targets the optimal MA risk, can be regarded as a variation of the sum of a set of positive-part Stein estimators from multiple mutually orthogonal subspaces. Specifically, both the MMA estimator and the positive-part Stein estimator share the common objective of minimizing unbiased risk estimation, albeit within different weight constraints. Second, we introduce a novel MA procedure to achieve asymptotic optimality by adapting the blockwise Stein rules from prior works \citep[][]{Donoho1995adapting, Nemirovski1998, Cavalier2001penalized} to linear regression. In particular, when the candidate model set is properly constructed, this Stein-type MA estimator achieves the full potential of MA (i.e., the minimal MA risk over all the nested models) in a sufficiently large parameter space. The results of finite sample simulations support our theories.

This paper gives the opportunity of looking at MA from different angles. By connecting MA with shrinkage, existing knowledge and technology derived from shrinkage estimation can be potentially transferred to MA. The selected review presented in this paper and the unveiled connections provide a theoretical foundation for this transfer; See Section~\ref{sec:discs} for more discussion.

The remainder of the paper is structured as follows. In Section~\ref{sec:setup}, we set up our regression problem. Section~\ref{sec:connect} draws the theoretical connections between MA and shrinkage. In Section~\ref{sec:stein}, we propose a Stein-type MA procedure and present its theoretical properties. Section~\ref{sec:simu} examines the finite sample properties of proposed methods by numerical simulations, and Section~\ref{sec:discs} concludes the paper. Proofs are included in the Appendix.

\section{Problem setup}\label{sec:setup}

Consider the linear regression model
\begin{equation}\label{eq:model}
  y_i =\mu_i+ \varepsilon_i=\sum_{j=1}^{p_n}\beta_jx_{ij} + \varepsilon_i, \quad i=1,\ldots,n,
\end{equation}
where $\varepsilon_1,\ldots,\varepsilon_n$ are i.i.d. Gaussian random errors with mean $0$ and variance $\sigma^2$, $x_{i1},\ldots,x_{ip_n}$, $i=1,\ldots,n$, are non-stochastic regressors, and $p_n$ is the number of regressors. Define $\by=(y_1,\ldots,y_n)^{\top}$, $\bmu=(\mu_1,\ldots,\mu_n)^{\top}$, $\boldsymbol{\varepsilon}=(\varepsilon_1,\ldots,\varepsilon_n)^{\top}$, $\bbeta = (\beta_1,\ldots,\beta_{p_n})^{\top}$, and $\bx_j = (x_{1j},\ldots,x_{nj})^{\top}$ for $j=1,\ldots,p_n$. Let $\bX=\left[\bx_1,\ldots,\bx_{p_n}\right]$ denote the regressor matrix. We can write (\ref{eq:model}) in matrix form
\begin{equation}\label{eq:model_m}
  \by =\bmu+ \boldsymbol{\varepsilon}= \bX \bbeta + \boldsymbol{\varepsilon}.
\end{equation}
Throughout this paper, we assume that $p_n\leq n$, and the regressor matrix $\bX$ has full column rank. For simplicity, we consider the case where the error variance $\sigma^2$ is known.

To estimate the true regression mean vector $\bmu$, $M_n$ strictly nested linear models are considered as candidates. The $m$-th candidate model includes the first $k_m$ regressors, where $1\leq k_1<k_2<\cdots<k_{M_n}\leq p_n$. The information about the sizes of candidate models is stored in a set $\mathcal{M}=\left\{k_1,k_2,\ldots,k_{M_n} \right\}$ with $M_n=|\mathcal{M}|$, where $|\mathcal{S}|$ denotes the cardinality of a set $\mathcal{S}$. Let $\bX_{k_m}=\left[\bx_1,\ldots,\bx_{k_m}\right]$ be the design matrix of the $m$-th candidate model. We estimate the regression coefficients by the least-squares method. The $m$-th estimator of $\bmu$ is
$$
\hat{\bmu}_{k_m}=\bX_{k_m}(\bX_{k_m}^{\top}\bX_{k_m})^{-1}\bX_{k_m}^{\top}\by = \bP_{k_m}\by,
$$
where $\bP_{k_m}\triangleq \bX_{k_m}(\bX_{k_m}^{\top}\bX_{k_m})^{-1}\bX_{k_m}^{\top}$ and $\bP_{0}\triangleq\boldsymbol{0}_{n \times n}$.

Let $\bw=(w_1,\ldots,w_{M_n})^{\top}$ denote a weight vector in the unit simplex of $\mathbb{R}^{M_n}$:
\begin{equation}\label{eq:weight_general}
  \mathcal{W}_{M_n}\triangleq\left\{\bw\in[0,1]^{M_n}:\sum_{m=1}^{M_n}w_m=1\right\}.
\end{equation}
Given the candidate model set $\mathcal{M}$, the MA estimator of $\bmu$ is
\begin{equation}\label{eq:MA}
  \hat{\bmu}_{\bw | \mathcal{M}}\triangleq\sum_{m=1}^{M_n}w_m\hat{\bmu}_{k_m}=\sum_{m=1}^{M_n}w_m\bP_{k_m}\by,
\end{equation}
where the subscript $\bw | \mathcal{M}$ is to emphasize the dependence of the MA estimator on the candidate model set $\mathcal{M}$.

To measure the performance of an estimator $\hat{\bmu}$, we consider the squared $\ell_2$ loss $L_n(\hat{\bmu}, \bmu)= \| \hat{\bmu} - \bmu \|^2$ and its corresponding risk $R_n(\hat{\bmu}, \bmu)=\mathbb{E}L_n(\hat{\bmu}, \bmu)$. For abbreviation, let $L_n(\bw | \mathcal{M}, \bmu)$ and $R_n(\bw | \mathcal{M}, \bmu)$ denote $L_n(\hat{\bmu}_{\bw | \mathcal{M}},\bmu)$ and $R_n(\hat{\bmu}_{\bw | \mathcal{M}},\bmu)$, respectively. The optimal weight vector in $\mathcal{W}_{M_n}$ is defined by
\begin{equation}\label{eq:optimal_weight_standard}
  \bw^*|\mathcal{M}\triangleq\argmin_{\bw \in \mathcal{W}_{M_n}}R_n(\bw | \mathcal{M}, \bmu).
\end{equation}

In this paper, our objective is to construct MA estimators which perform asymptotically equivalent to (or even better than) that based on $\bw^*|\mathcal{M}$. The specific definitions are as follows.
\begin{definition}\label{def:aop}
  An MA estimator with the weights $\hat{\bw}|\mathcal{M}$ trained on data is called asymptotically optimal if
  \begin{equation}\label{eq:aop}
  \mathbb{E}L_n(\hat{\bw} | \mathcal{M}, \bmu) \leq \left[ 1+o(1) \right]R_n(\bw^* | \mathcal{M}, \bmu)
\end{equation}
holds as $n \to \infty$.
\end{definition}
A representative example of the candidate model set is $\mathcal{M}_a=\{1,2,\ldots,p_n \}$, which encompasses all nested models. For any $\mathcal{M}\subseteq \mathcal{M}_a$, it holds that
$$
R_n(\bw^* | \mathcal{M}_a,\bmu)\leq R_n(\bw^* | \mathcal{M},\bmu).
$$
Thus, $R_n(\bw^* | \mathcal{M}_a,\bmu)$ represents the full potential of MA in our context.
\begin{definition}\label{def:fullaop}
An MA estimator with the weights $\hat{\bw}|\mathcal{M}$ estimated on data is called fully asymptotically optimal if
\begin{equation}\label{eq:fullaop}
  \mathbb{E}L_n(\hat{\bw} | \mathcal{M}, \bmu) \leq \left[ 1+o(1) \right]R_n(\bw^* | \mathcal{M}_a, \bmu)
\end{equation}
holds as $n \to \infty$.
\end{definition}
Note that Definition~\ref{def:fullaop} focuses on the full potential of MA in the nested model setting we considered, while Definition~\ref{def:aop} provides the asymptotic justification for the MA estimators based on general nested model sets $\mathcal{M}$. In the next sections, we demonstrate some unexplored connections between MA and shrinkage estimation. Through these connections, we develop several Stein-type MA procedures to achieve the asymptotic optimality properties outlined in Definitions~\ref{def:aop}--\ref{def:fullaop}.

\section{Connecting MA to shrinkage}\label{sec:connect}

\subsection{Optimal MA and monotone oracle}\label{subsec:MA_monotone}

We first introduce some notations that will play a key role in our theoretical analysis. Given a nested candidate model set $\mathcal{M}=\left\{k_1,k_2,\ldots,k_{M_n} \right\}$, define the matrixes $\bD_{m|\mathcal{M}}\triangleq\bP_{k_m}-\bP_{k_{m-1}}$ for $m=1,\ldots,M_n$, and $\bD_{M_n+1|\mathcal{M}}\triangleq\bI - \bP_{k_{M_n}}$. As pointed out in \cite{Xu2022From}, $\bD_{m|\mathcal{M}},m=1,\ldots,M_n+1$ are mutually orthogonal since $\bD_{m|\mathcal{M}}\bD_{m'|\mathcal{M}}=\bD_{m|\mathcal{M}}\delta_{mm'}$, where $\delta_{mm'}$ is the Kronecker delta.

We represent the response vector $\by$ in $M_n+1$ orthogonal subspaces
\begin{equation}\label{eq:model_s}
  \by_{m|\mathcal{M}} =\bmu_{m|\mathcal{M}}+ \boldsymbol{\varepsilon}_{m|\mathcal{M}}, \quad m=1,\ldots,M_n+1,
\end{equation}
where $\by_{m|\mathcal{M}}\triangleq\bD_{m|\mathcal{M}}\by$, $\bmu_{m|\mathcal{M}} \triangleq \bD_{m|\mathcal{M}}\bmu$, and $\boldsymbol{\varepsilon}_{m|\mathcal{M}} \triangleq \bD_{m|\mathcal{M}}\boldsymbol{\varepsilon}$. The expression (\ref{eq:model_s}) defines a Gaussian sequence model with $M_n+1$ independent vector-valued observations. Note that $\boldsymbol{\varepsilon}_{m|\mathcal{M}}$ is distributed as $N(\boldsymbol{0}, \sigma^2\bD_{m|\mathcal{M}} )$ and is independent of $\boldsymbol{\varepsilon}_{m'|\mathcal{M}}$ when $m'\neq m$. The vector $\bmu_{m|\mathcal{M}}$ can be interpreted as the regression mean vector explained solely by the $m$-th candidate model among the first $m$ nested models.

Based on the Gaussian sequence model (\ref{eq:model_s}), the MA estimator (\ref{eq:MA}) can be written as a linear combination of $\by_{m|\mathcal{M}}$
\begin{equation}\label{eq:MA_ga}
  \hat{\bmu}_{\bw | \mathcal{M}}=\hat{\bmu}_{\bgamma| \mathcal{M}}\triangleq \sum_{m=1}^{M_n} \gamma_m \by_{m|\mathcal{M}},
\end{equation}
where $\gamma_m = \sum_{j=m}^{M_n}w_j$, $\gamma_{M_n+1}=0$, and $\bgamma|\mathcal{M}=(\gamma_1,\ldots,\gamma_{M_n})^{\top}$ is the cumulative weight vector. In a slight abuse of notation, we simultaneously write $\bw | \mathcal{M}$ and $\bgamma| \mathcal{M}$ in the subscript of $\hat{\bmu}$ in (\ref{eq:MA_ga}). It will cause no confusion if we use $\bw | \mathcal{M}$ and $\bgamma| \mathcal{M}$ to designate the forms of MA in terms of the weights and the cumulative weights, respectively. Similar arguments apply to the other notations defined in this section.

We observe that the MA risk equals to the risk of the linear estimator (\ref{eq:MA_ga}) in the Gaussian sequence model (\ref{eq:model_s}). That is
\begin{equation}\label{eq:risk}
\begin{split}
   R_n(\bw | \mathcal{M}, \bmu) &=R_n(\bgamma | \mathcal{M}, \bmu)\\
   &\triangleq \sum_{m=1}^{M_n+1}\mathbb{E}\| \bmu_{m|\mathcal{M}} - \gamma_m \by_{m|\mathcal{M}} \|^2 \\
     & = \sum_{m=1}^{M_n}\left[ \|\bmu_{m|\mathcal{M}}\|^2 ( 1 - \gamma_m)^2  + \sigma^2_{m|\mathcal{M}}\gamma_m^2 \right]+ \|\bmu_{M_n+1|\mathcal{M}}\|^2,
\end{split}
\end{equation}
where $\sigma^2_{m|\mathcal{M}}\triangleq(k_m-k_{m-1})\sigma^2$. Additionally, restricting the weight vector in $\mathcal{W}_{M_n}$ is equivalent to requiring the cumulative weight vector to belong to the set
\begin{equation}\label{eq:monotone}
  \Gamma_{M_n} \triangleq\left\{\bgamma :1=\gamma_1\geq \gamma_2 \geq \cdots \geq \gamma_{M_n}\geq 0\right\}.
\end{equation}

We define the monotone oracle as
$$
\bgamma^*|\mathcal{M}\triangleq\argmin_{\bgamma \in \Gamma_{M_n}}R_n(\bgamma | \mathcal{M}, \bmu).
$$
Therefore, the optimal weight vector $\bw^*|\mathcal{M}=(w_1^*,\ldots,w_{M_n}^*)^{\top}$ connects to $\bgamma^*|\mathcal{M}$ by $w_m^*=\gamma_m^*-\gamma_{m+1}^*$ for $m=1,\ldots,M_n$ and $\gamma_{M_n+1}^*=0$.

\begin{remark}
 It is noteworthy that the connections between the optimal MA estimator and the monotone oracle $\bgamma^*|\mathcal{M}$ have also been utilized in \cite{Peng2022improvability} and \cite{Xu2022From}. However, these two papers mainly investigate the property of (\ref{eq:risk}) rather than focusing on specific MA procedures. In the subsequent subsection, we will investigate how to estimate the optimal MA estimator or the monotone oracle based on the observed data, providing further insight into the connection between shrinkage and MA.
\end{remark}

\subsection{MMA and Stein estimation}\label{sec:MMA_stein}

A well-established approach for estimating $\bw^*|\mathcal{M}$ and $\bgamma^*|\mathcal{M}$ is based on the principle of unbiased risk estimation (URE), which has been deeply studied in the literature \citep[see, e.g.,][]{Mallows1973, Akaike1973, Stein1981Estimation, Golubev1990, Beran1998modulation}. Under the linear regression model (\ref{eq:model_m}), the principle of URE is known as the MMA criterion
\begin{equation}\label{eq:criterion}
  C_n(\bw|\mathcal{M}, \by)=\| \by - \hat{\bmu}_{\bw | \mathcal{M}} \|^2 + 2\sigma^2\bk^{\top}\bw.
\end{equation}
Minimizing $C_n(\bw|\mathcal{M}, \by)$ over $\mathcal{W}_{M_n}$ yields $\hat{\bw}_{\text{\scriptsize MMA} }|\mathcal{M}=\arg\min_{\bw \in \mathcal{W}_{M_n}}C_n(\bw|\mathcal{M}, \by)$ and the MMA estimator
\begin{equation}\label{eq:class_MMA}
  \hat{\bmu}_{\hat{\bw}_{\text{\tiny MMA} }|\mathcal{M}}=\sum_{m=1}^{M_n}\hat{w}_{\text{\tiny MMA},m }\hat{\bmu}_{k_m},
\end{equation}
where $\hat{w}_{\text{\tiny MMA},m }$ denotes the $m$-th element of $\hat{\bw}_{\text{\scriptsize MMA} }|\mathcal{M}$.

Under the Gaussian sequence model (\ref{eq:model_s}), an equivalent expression of the MMA criterion in terms of $\bgamma|\mathcal{M}$ is
\begin{equation}\label{eq:criterion2}
  C_n(\bgamma|\mathcal{M}, \by)=\sum_{m=1}^{M_n}\left[ \|\by_{m|\mathcal{M}}\|^2 (1-\gamma_m)^2  + 2\sigma^2_{m|\mathcal{M}}\gamma_m  \right]+ \|\by_{M_n+1|\mathcal{M}} \|^2,
\end{equation}
where $\by_{m|\mathcal{M}}$ and $\sigma^2_{m|\mathcal{M}}$ are defined in Section~\ref{subsec:MA_monotone}. Minimizing $C_n(\bgamma|\mathcal{M}, \by)$ over the monotone weight set $\Gamma_{M_n}$ gives $\hat{\bgamma}_{\text{\scriptsize MMA} }|\mathcal{M}=\min_{\bgamma \in \Gamma_{M_n}}C_n(\bgamma|\mathcal{M}, \by)$ and
$$
\hat{\bmu}_{\hat{\bgamma}_{\text{\tiny MMA} }|\mathcal{M}}=\sum_{m=1}^{M_n} \hat{\gamma}_{\text{\tiny MMA},m } \by_{m|\mathcal{M}},
$$
where $\hat{\gamma}_{\text{\scriptsize MMA},m }$ is the $m$-th element of $\hat{\bgamma}_{\text{\scriptsize MMA} }|\mathcal{M}$.
Based on the relation (\ref{eq:MA_ga}), it is evident that $\hat{\bmu}_{\hat{\bw}_{\text{\tiny MMA} }|\mathcal{M}}$ and $\hat{\bmu}_{\hat{\bgamma}_{\text{\tiny MMA}}|\mathcal{M}}$ are exactly the same estimator.

In general, $\hat{\bmu}_{\hat{\bw}_{\text{\tiny MMA} }|\mathcal{M}}$ and $\hat{\bmu}_{\hat{\bgamma}_{\text{\tiny MMA} }|\mathcal{M}}$ do not possess explicit expressions due to the positive and equality constraints in  $\mathcal{W}_{M_n}$ and the monotonicity constraints in $\Gamma_{M_n}$, respectively. However, if we consider minimizing $C_n(\bgamma|\mathcal{M}, \by)$ within a hypercube
$$
\tilde{\Gamma}_{M_n}\triangleq[0,1]^{M_n}
$$
instead of the monotone set $\Gamma_{M_n}$, or equivalently, minimizing $C_n(\bw|\mathcal{M}, \by)$ over a relaxed weight set
\begin{equation}\label{eq:W_tilde}
  \tilde{\mathcal{W}}_{M_n} \triangleq \left\{ \bw \in \mathbb{R}^{M_n}:0\leq \sum_{j=m}^{M_n}w_j \leq 1,m=1,\ldots,M_n\right\}
\end{equation}
rather than $\mathcal{W}_{M_n}$, then the problem becomes more tractable. The solution to the optimization problem $\min_{\bgamma \in \tilde{\Gamma}_{M_n}}C_n(\bgamma|\mathcal{M}, \by)$ takes an explicit form denoted as  $\hat{\bgamma}_{\text{\scriptsize STE} }|\mathcal{M}=(\hat{\gamma}_{\text{\scriptsize STE},1},\ldots, \hat{\gamma}_{\text{\scriptsize STE},M_n})^{\top}$, where
\begin{equation}\label{eq:gamma_1}
  \hat{\gamma}_{\text{\scriptsize STE},m} \triangleq \left(1 - \frac{\sigma^2_{m|\mathcal{M}}}{\|\by_{m|\mathcal{M}}\|^2}\right)_{+}, \quad m=1,\ldots, M_n,
\end{equation}
and $x_+$ denotes $\max(x,0)$.

Define the weight vector $\hat{\bw}_{\text{\scriptsize STE} }|\mathcal{M}=(\hat{w}_{\text{\scriptsize STE},1},\ldots,\hat{w}_{\text{\scriptsize STE},M_n})^{\top}$, where $\hat{w}_{\text{\scriptsize STE},m}=\hat{\gamma}_{\text{\scriptsize STE},m}-\hat{\gamma}_{\text{\scriptsize STE},m+1}$ for $m=1,\ldots,M_n$, and $\hat{\gamma}_{\text{\scriptsize STE},M_n+1}=0$. The corresponding MA estimator is given by
\begin{equation}\label{eq:est_shr}
\begin{split}
    \hat{\bmu}_{\hat{\bw}_{\text{\tiny STE} }|\mathcal{M}}
&= \hat{\bmu}_{\hat{\bgamma}_{\text{\tiny STE} }|\mathcal{M}} = \sum_{m=1}^{M_n} \left(1 - \frac{\sigma^2_{m|\mathcal{M}}}{\|\by_{m|\mathcal{M}}\|^2}\right)_{+} \by_{m|\mathcal{M}},
\end{split}
\end{equation}
which actually is the sum of multiple positive-part Stein estimators in different orthogonal subspaces generated by $\bD_{m|\mathcal{M}},m=1,\ldots,M_n$.

  Comparing the MMA estimator $\hat{\bmu}_{\hat{\bw}_{\text{\tiny MMA} }|\mathcal{M}}$ with the Stein estimator $\hat{\bmu}_{\hat{\bw}_{\text{\tiny STE} }|\mathcal{M}}$, we observe that both serve as the minimizers of the principle of URE, albeit under different weight constraints. Specifically, $\hat{\bmu}_{\hat{\bw}_{\text{\tiny MMA} }|\mathcal{M}}$ minimizes the principle of URE within the unit simplex $\mathcal{W}_{M_n}$, while $\hat{\bmu}_{\hat{\bw}_{\text{\tiny STE} }|\mathcal{M}}$ is based on the relaxed weight set $\tilde{\mathcal{W}}_{M_n}$. The weight relaxation from $\mathcal{W}_{M_n}$ to $\tilde{\mathcal{W}}_{M_n}$ can yield significant benefits for the nested MA under certain conditions (see Section~\ref{subsec:relax}) while incurring controllable cost to achieve asymptotic optimality (see Section~\ref{sec:stein}). %It is worthy mentioning that although there is no general closed form for $\hat{\bmu}_{\hat{\bw}_{\text{\tiny MMA} }|\mathcal{M}}$ in the multiple model context, investigating its properties directly is still feasible \citep[see, e.g.,][]{Hansen2014risk, Zhang2016dominance, Zhang_liu_2019, Peng2023optimality}.

\begin{remark}

In a nested model setting similar to ours, \cite{Hansen2014risk} reformulated the MA estimator and the MMA criterion in term of the cumulative weight vector $\blambda = (\lambda_1,\ldots,\lambda_{M_n})^{\top}$, where $\lambda_m = \sum_{j=1}^{m}w_j,m=1,\ldots,M_n$. In contrast, we rewrite the nested MA estimator in terms of $\bgamma = (\gamma_1,\ldots,\gamma_{M_n})^{\top}$ with $\gamma_m = \sum_{j=m}^{M_n}w_j, m=1,\ldots,M_n$. Notably, given the same candidate model set and weight constraint, the formulas of MA expressed in terms of $\blambda$ and $\bgamma$ are mathematically equivalent due to the one-to-one correspondence
\begin{equation}\label{eq:one_to_one}
  \blambda=\begin{pmatrix}
             1 & -1 & 0 &  \cdots & 0 \\
             1 &  0  & -1 & \cdots  & 0 \\
             \vdots  & \vdots  &  \vdots &   &  \vdots \\
             1 & 0 & 0 &  \cdots & -1\\
             1 & 0 & 0 &  \cdots & 0
           \end{pmatrix}\bgamma.
\end{equation}
With the monotonically non-decreasing constrains $0 \leq \lambda_1 \leq \cdots \leq \lambda_{M_n}=1$ as imposed in \cite{Hansen2014risk}, the closed forms for the minimizers of MMA are derived only in some special cases, such as the two-model scenario considered in Section~6 of \cite{Hansen2014risk}. Based on the relation (\ref{eq:one_to_one}) and some additional algebra operations, it is evident that minimizing \cite{Hansen2014risk}'s MMA formula over the relaxed constraints $0 \leq \lambda_{M_n} \leq 1$ and $0 \leq \lambda_{M_n} - \lambda_{m} \leq 1, 1\leq m \leq M_n-1$ results in the same positive-part Stein estimator (\ref{eq:est_shr}).

\end{remark}

\subsection{Relaxation of weight constraint}\label{subsec:relax}

The results in the last subsection reveal that the standard MMA estimator is different from the Stein estimator (\ref{eq:est_shr}) in weight constraint. An intriguing question arises: Does relaxing the weight constraint from the simplex $\mathcal{W}_{M_n}$ to $\tilde{\mathcal{W}}_{M_n}$ yield substantial benefits for MA? This subsection addresses this question by comparing the optimal MA risks within $\mathcal{W}_{M_n}$ and $\tilde{\mathcal{W}}_{M_n}$.

Recall the definition in (\ref{eq:optimal_weight_standard}) that $\bw^*|\mathcal{M}$ denotes the optimal weight vector within $\mathcal{W}_{M_n}$. Additionally, define
\begin{equation*}
  \tilde{\bw}^*|\mathcal{M}\triangleq\argmin_{\bw \in \tilde{\mathcal{W}}_{M_n}}R_n(\bw | \mathcal{M}, \bmu)
\end{equation*}
as the optimal weight vector in $\tilde{\mathcal{W}}_{M_n}$. Since $\mathcal{W}_{M_n}\subseteq \tilde{\mathcal{W}}_{M_n}$, the relation
\begin{equation*}
  R_n(\tilde{\bw}^*|\mathcal{M}, \bmu)\leq R_n(\bw^*|\mathcal{M}, \bmu)
\end{equation*}
always holds for any candidate model set $\mathcal{M}$ and regression mean vector $\bmu$.

To further elucidate the difference between $R_n(\tilde{\bw}^*|\mathcal{M}, \bmu)$ and $R_n(\bw^*|\mathcal{M}, \bmu)$, define
\begin{equation}\label{eq:M_T_def}
  \mathcal{M}_T\triangleq \argmax_{\{k_1,k_{M_n} \}\subseteq \mathcal{S}\subseteq \mathcal{M}}\left\{ \frac{ \|\bmu_{l|\mathcal{S}}\|^2}{\sigma^2_{l|\mathcal{S}}}\geq  \frac{ \|\bmu_{l+1|\mathcal{S}}\|^2}{\sigma^2_{l+1|\mathcal{S}}},\, l=2,\ldots, |\mathcal{S}|-1\right\}
\end{equation}
if such a set exists. Otherwise, let $\mathcal{M}_T$ be $\{k_1,k_{M_n} \}$. Define $L_n = |\mathcal{M}_T|$ as the number of candidate models contained in $\mathcal{M}_T$. Based on the definition of $\mathcal{M}_T$, it is evident that
\begin{equation*}
  \{k_1,k_{M_n} \} \subseteq \mathcal{M}_T \subseteq \{k_1,k_2,\ldots,k_{M_n} \}
\end{equation*}
and $2 \leq L_n \leq M_n$. %Section~\ref{sec:proof:improv} provides an iterative approach to identify the set $\mathcal{M}_T$.
It is worth noting that $\|\bmu_{l|\mathcal{S}}\|^2/\sigma^2_{l|\mathcal{S}}$ represents the signal-to-noise ratio (SNR) in the $l$-th subspace generated by the candidate model set $\mathcal{S}$. Thus, $\mathcal{M}_T$ is the largest subset of $\mathcal{M}$ whose SNRs from the second to last subspaces are monotonically nonincreasing.

\begin{theorem}\label{theo:improv}
Given a candidate model set $\mathcal{M}=\{k_1,k_2,\ldots,k_{M_n} \}$, we have
\begin{equation}\label{eq:optimal_enlarge}
  R_n(\tilde{\bw}^*|\mathcal{M}, \bmu)=\sum_{m=1}^{M_n}\frac{\|\bmu_{m|\mathcal{M}}\|^2\sigma^2_{m|\mathcal{M}}}{\|\bmu_{m|\mathcal{M}}\|^2+\sigma^2_{m|\mathcal{M}}}+ \|\bmu_{M_n+1|\mathcal{M}}\|^2,
\end{equation}
and
\begin{equation}\label{eq:optimal_convex}
\begin{split}
    R_n(\bw^*|\mathcal{M}, \bmu)&=R_n(\bw^*|\mathcal{M}_T, \bmu) \\
     & =\sigma^2_{1|\mathcal{M}_T}+\sum_{l=2}^{L_n}\frac{\|\bmu_{l|\mathcal{M}_T}\|^2\sigma^2_{l|\mathcal{M}_T}}{\|\bmu_{l|\mathcal{M}_T}\|^2+\sigma^2_{l|\mathcal{M}_T}}+ \|\bmu_{L_n+1|\mathcal{M}_T}\|^2.
\end{split}
\end{equation}

\end{theorem}

Theorem~\ref{theo:improv} gives explicit expressions for the optimal MA risks $R_n(\tilde{\bw}^*|\mathcal{M}, \bmu)$ and $R_n(\bw^*|\mathcal{M}, \bmu)$. If $\|\bmu_{l|\mathcal{M}}\|^2/\sigma^2_{l|\mathcal{M}}\geq \|\bmu_{l+1|\mathcal{M}}\|^2/\sigma^2_{l+1|\mathcal{M}}$ for all $l=2,\ldots,M_n-1$, then $\mathcal{M}_T=\mathcal{M}=\{k_1,k_2,\ldots,k_{M_n} \}$, and
\begin{equation*}
  R_n(\bw^*|\mathcal{M}, \bmu)-R_n(\tilde{\bw}^*|\mathcal{M}, \bmu)=\frac{\sigma^4_{1|\mathcal{M}}}{\|\bmu_{1|\mathcal{M}}\|^2+\sigma^2_{1|\mathcal{M}}},
\end{equation*}
which is negligible compared to $R_n(\bw^*|\mathcal{M}, \bmu)$ provided $k_1$ is bounded and $R_n(\bw^*|\mathcal{M}, \bmu) \to \infty$. This result implies that if the SNR $\|\bmu_{l|\mathcal{M}}\|^2/\sigma^2_{l|\mathcal{M}}$ is monotonically nonincreasing as $l$ increases, the optimal MA risks in $\mathcal{W}_{M_n}$ and $\tilde{\mathcal{W}}_{M_n}$ are asymptotically equivalent.

However, when the SNRs in the subspaces generated by the original model set $\mathcal{M}$ are not monotonically non-increasing, the weight restriction with $\mathcal{W}_{M_n}$ may limit the potential of MA since $R_n(\bw^*|\mathcal{M}, \bmu)$ equals to the optimal MA risk based on the reduced candidate models set $\mathcal{M}_T$. To compare the rate of $R_n(\tilde{\bw}^*|\mathcal{M}, \bmu)$ with the rate of $R_n(\bw^*|\mathcal{M}_T, \bmu)$, we introduce several additional notations. Let $\|\bmu_{(1)|\mathcal{M}}\|^2/\sigma^2_{(1)|\mathcal{M}}, \ldots, \|\bmu_{(M_n)|\mathcal{M}}\|^2/\sigma^2_{(M_n)|\mathcal{M}}$ be a non-increasing rearrangement of $\|\bmu_{1|\mathcal{M}}\|^2/\sigma^2_{1|\mathcal{M}}, \ldots, \|\bmu_{M_n|\mathcal{M}}\|^2/\sigma^2_{M_n|\mathcal{M}}$. Define
\begin{equation*}
  m_n^* \triangleq \left|\left\{m: \|\bmu_{(m)|\mathcal{M}}\|^2/\sigma^2_{(m)|\mathcal{M}} \geq  1 \right\} \right|.
\end{equation*}
Since $\|\bmu_{2|\mathcal{M}_T}\|^2/\sigma^2_{2|\mathcal{M}_T}, \ldots, \|\bmu_{L_n|\mathcal{M}_T}\|^2/\sigma^2_{L_n|\mathcal{M}_T}$ are monotonically non-increasing, we define
\begin{equation*}
  l_n^* \triangleq \left|\left\{2 \leq l \leq L_n: \|\bmu_{l|\mathcal{M}_T}\|^2/\sigma^2_{l|\mathcal{M}_T} \geq 1 \right\} \right|+1.
\end{equation*}

\begin{corollary}\label{theo:improv2}
Given a candidate model set $\mathcal{M}=\{k_1,k_2,\ldots,k_{M_n} \}$ with $k_{M_n}=p_n$. If the regression mean vector $\bmu$ satisfies
\begin{equation}\label{eq:compare_cond}
  \sum_{m=m_n^*+1}^{M_n}\|\bmu_{(m)|\mathcal{M}}\|^2=O\left(\sum_{m=1}^{m_n^*}\sigma^2_{(m)|\mathcal{M}}\right),
\end{equation}
then we have
  \begin{equation}\label{eq:upper_bound_ratio}
    \frac{R_n(\tilde{\bw}^*|\mathcal{M}, \bmu)}{R_n(\bw^*|\mathcal{M}, \bmu)} = O\left( \frac{\sum_{m=1}^{m_n^*}\sigma^2_{(m)|\mathcal{M}}}{\sum_{l=1}^{l_n^*}\sigma^2_{l|\mathcal{M}_T}}\right),
  \end{equation}
  where $\sigma^2_{(m)|\mathcal{M}}$ is the variance component in $\|\bmu_{(m)|\mathcal{M}}\|^2/\sigma^2_{(m)|\mathcal{M}}$.
\end{corollary}
Corollary~\ref{theo:improv2} provides a general upper bound on the rate of ratio between $R_n(\tilde{\bw}^*|\mathcal{M}, \bmu)$ and $R_n(\bw^*|\mathcal{M}, \bmu)$. To illustrate this corollary, assume $p_n=n$ and consider the candidate model set with equally increasing size $k_m = m k_1 $ for $m = 1,\ldots,M_n-1$ and $k_{M_n} = n$, where $k_1$ is a fixed integer, and $M_n=\argmin_{m}\{ m k_1 > n \}$. Then, based on the first inequality in the proof of Theorem~2 in \cite{Peng2022improvability}, the regularity condition (\ref{eq:compare_cond}) can be satisfied
%if there exists a constant $\kappa>1$ and $0<\nu < 1$ with $\kappa \nu^2 <1$ such that
%\begin{equation}\label{eq:compare_2}
%  \frac{\|\bmu_{(\lfloor \kappa l \rfloor)|\mathcal{M}}\|^2/\sigma^2_{(\lfloor \kappa l \rfloor)|\mathcal{M}}}{\|\bmu_{( l )|\mathcal{M}}\|^2/\sigma^2_{( l )|\mathcal{M}}} \leq \nu
%\end{equation}
%when $l$ is large enough. Actually, the condition (\ref{eq:compare_2}) can be satisfied
when $\|\bmu_{( l )|\mathcal{M}}\|^2/\sigma^2_{( l )|\mathcal{M}}$ decays at a polynomial rate
\begin{equation*}
  \frac{\|\bmu_{( l )|\mathcal{M}}\|^2}{n\sigma^2_{( l )|\mathcal{M}}} \asymp l^{-2\alpha_1}, \, \frac{1}{2}< \alpha_1 < \infty
\end{equation*}
or at an exponential rate
\begin{equation*}
  \frac{\|\bmu_{( l )|\mathcal{M}}\|^2}{n\sigma^2_{( l )|\mathcal{M}}} \asymp \exp(-2l^{\alpha_2 } ), \, 0< \alpha_2 < \infty.
\end{equation*}

Consider a specific case where $\|\bmu_{l|\mathcal{M}}\|^2/\sigma^2_{l|\mathcal{M}}< \|\bmu_{l+1|\mathcal{M}}\|^2/\sigma^2_{l+1|\mathcal{M}}$ for every $l=2,\ldots,M_n-1$, and $\lim\inf_n \|\bmu_{2|\mathcal{M}_T}\|^2/\sigma^2_{2|\mathcal{M}_T} \geq 1$. Then $\mathcal{M}_T$ is reduced to $\{k_1,n \}$, and $l_n^*=2$. In this case, we see
\begin{equation*}
   \frac{\sum_{m=1}^{m_n^*}\sigma^2_{(m)|\mathcal{M}}}{\sum_{l=1}^{l_n^*}\sigma^2_{l|\mathcal{M}_T}} = \frac{m_n^*k_1\sigma^2}{n\sigma^2} \to 0
\end{equation*}
if $m_n^*=o(n)$, which is satisfied when $\|\bmu_{( l )|\mathcal{M}}\|^2/\sigma^2_{( l )|\mathcal{M}}$ decays in the both polynomial and exponential rates. Hence, in the cases where the monotonicity of SNRs is violated, the optimal MA risk under the relaxed weight constraint may converge significantly faster compared to that under the standard simplex constraint. Furthermore, $R_n(\tilde{\bw}^*|\mathcal{M}, \bmu)$ is indeed attained by the negative weights $\tilde{w}_m^*=\tilde{\gamma}_m^* - \tilde{\gamma}_{m+1}^*,m=1,\ldots,M_n$ in this scenario, since
\begin{equation*}
  \tilde{\gamma}_m^* = \frac{\|\bmu_{m|\mathcal{M}}\|^2}{\|\bmu_{m|\mathcal{M}}\|^2+ \sigma^2_{m|\mathcal{M}}}, \quad m=1,\ldots,M_n,
\end{equation*}
are monotonically increasing.

\section{A Stein-type MA procedure}\label{sec:stein}

\subsection{Penalized blockwise Stein method}\label{sec:pbs}

Based on the works of \cite{Cavalier2001penalized, Cavalier2002Sharp}, we extend the Stein rule (\ref{eq:gamma_1}) to a penalized blockwise Stein rule, formulated as
\begin{equation}\label{eq:stein2}
  \hat{\gamma}_m = \left(1 - \frac{\sigma_{m|\mathcal{M}}^2(1+\varphi_m)}{\|\by_{m|\mathcal{M}}\|^2}\right)_{+},\quad m=1,\ldots,M_n,
\end{equation}
where $0\leq\varphi_m<1$ is a penalty factor, $\hat{\bgamma}|\mathcal{M}=(\hat{\gamma}_1,\ldots,\hat{\gamma}_{M_n})^{\top}$, and $\hat{\gamma}_{M_n+1}=0$. Subsequently, we define a Stein-type MA estimator as:
\begin{equation}\label{eq:stein_ma}
  \hat{\bmu}_{\hat{\bw} | \mathcal{M}} = \sum_{m=1}^{M_n} \hat{\gamma}_m \by_{m|\mathcal{M}}=\sum_{m=1}^{M_n}\hat{w}_m\hat{\bmu}_{k_m},
\end{equation}
where the $m$-th element of $\hat{\bw}|\mathcal{M}$ is $\hat{w}_m=\hat{\gamma}_m-\hat{\gamma}_{m+1}$.

Notably, when $\varphi_m = 0$, the Stein-type MA estimator (\ref{eq:stein_ma}) reduces to (\ref{eq:est_shr}). In contrast, in cases where $\varphi_m > 0$, the inclusion of the penalty factor may lead to a reduction in the number of nonzero cumulative weights compared to the Stein estimator (\ref{eq:est_shr}), thereby leading to better adaptivity to the scenarios with sparse signals. In the majority of this section, $\varphi_m$ is presumed to be positive and correlated with $n$. We defer discussion of the specific case where $\varphi_m=0$ for $1\leq m\leq M_n$ until the end of Section~\ref{sec:pbs}.

Specifically, we make the assumption that the values of $\varphi_m$ are small and satisfy
\begin{equation}\label{eq:ass0}
  \max_{1\leq m \leq M_n}\varphi_m \to 0
\end{equation}
as $n \to \infty$. To derive the theoretical properties of the Stein-type MA estimator $\hat{\bmu}_{\hat{\bw} | \mathcal{M}}$, we require two additional assumptions concerning $\mathcal{M}$ and $\varphi_m$.

\begin{assumption}\label{ass:1}
There exists a constant $c_1$ such that
\begin{equation}
  \sum_{m=1}^{M_n}\exp\left[  -\frac{(k_m-k_{m-1})\varphi_m^2}{16(1+2\sqrt{\varphi_m})^2}\right]\leq c_1.
\end{equation}
\end{assumption}

\begin{assumption}\label{ass:2}
For all $m=1,\ldots, M_n$, assume
  \begin{equation}
  \frac{1}{k_m-k_{m-1}} \leq \frac{1-\varphi_m}{4}.
\end{equation}
\end{assumption}

A prerequisite for Assumption~\ref{ass:1} is
\begin{equation}\label{eq:pre}
  (k_m-k_{m-1})\varphi_m^2\to \infty.
\end{equation}
This necessitates that $k_m - k_{m-1} \to \infty$ as $m$ increase. Equation (\ref{eq:pre}) suggests a typical choice for the penalty factor as $\varphi_m=1/(k_m-k_{m-1})^{\tau}$, where $0<\tau<1/2$. As for Assumption~\ref{ass:2}, a sufficient condition is $k_m-k_{m-1} > 3$, which aligns with a common assumption in Stein-type methods \citep[see, e.g.,][]{Stein1981Estimation}.

\begin{theorem}\label{theorem:1}
  Suppose Assumptions~\ref{ass:1}--\ref{ass:2} hold. Then for any sample size $n$ and any regression mean vector $\bmu$, we have
  \begin{equation}\label{eq:bound_penal_stein}
  \mathbb{E}L_n(\hat{\bw} | \mathcal{M},\bmu) \leq   (1+\bar{\varphi}) R_n(\tilde{\bw}^* | \mathcal{M},\bmu)+8c_1\sigma^2,
\end{equation}
where $\bar{\varphi} = \max_{1\leq m \leq M_n}\{ 2\varphi_m+16/[(k_m-k_{m-1})\varphi_m]\}$ and $R_n(\tilde{\bw}^* | \mathcal{M},\bmu)=\min_{\bw \in \tilde{\mathcal{W}}_{M_n}}R_n(\bw | \mathcal{M}, \bmu)$ denotes the optimal MA risk in the relaxed weight set $\tilde{\mathcal{W}}_{M_n}$.
\end{theorem}

The oracle inequality provided in Theorem~\ref{theorem:1} implies that the Stein-type MA estimator (\ref{eq:stein_ma}) asymptotically performs as well as the optimal MA estimator based on $\mathcal{M}$ and $\tilde{\mathcal{W}}_{M_n}$. Specifically, when both (\ref{eq:ass0}) and (\ref{eq:pre}) hold,
\begin{equation}\label{eq:bar_phi}
  \bar{\varphi} = \max_{1\leq m \leq M_n}\left[ 2\varphi_m+o(\varphi_m)\right]=O\left(\max_{1\leq m \leq M_n}\varphi_m\right)\to 0.
\end{equation}
Thus, if $R_n(\tilde{\bw}^* | \mathcal{M},\bmu)\to \infty$ as $n \to \infty$, (\ref{eq:stein_ma}) is asymptotically optimal in the sense that
\begin{equation*}
  \mathbb{E}L_n(\hat{\bw} | \mathcal{M},\bmu) \leq [1+o(1)] R_n(\tilde{\bw}^* | \mathcal{M},\bmu).
\end{equation*}
Moreover, as $\mathcal{W}_{M_n}\subseteq\tilde{\mathcal{W}}_{M_n}$ and $R_n(\tilde{\bw}^* | \mathcal{M},\bmu) \leq R_n(\bw^* | \mathcal{M},\bmu)$, Theorem~\ref{theorem:1} also establishes the asymptotic optimality of (\ref{eq:stein_ma}) in terms of (\ref{eq:aop}).

Note that Assumption~\ref{ass:1} excludes the case where $\varphi_m = 0, 1 \leq m \leq M_n$.  Consequently, the risk bound presented in Theorem~\ref{theorem:1} does not apply to the non-penalized blockwise Stein estimator (\ref{eq:est_shr}). However, by making slight modifications to the proof technique in Theorem~\ref{theorem:1} (see Section~\ref{sec:proof_cor_41} of the Appendix), we derive the following bound on the risk of (\ref{eq:est_shr}).

\begin{assumption}\label{ass:2_add}
For all $m=1,\ldots, M_n$, assume $k_m-k_{m-1} > 3$.
\end{assumption}

\begin{corollary}\label{cor:classical_stein}
  Suppose Assumption~\ref{ass:2_add} is satisfied, then we have
  \begin{equation}\label{eq:bound_class_stein}
    \mathbb{E}L_n(\hat{\bw}_{\text{\scriptsize STE} }|\mathcal{M},\bmu) \leq R_n(\tilde{\bw}^*|\mathcal{M}, \bmu) + 4M_n\sigma^2.
  \end{equation}
\end{corollary}

Assumption~\ref{ass:2_add} serves as a foundational assumption for the classical James-Stein estimator \citep{Stein1981Estimation}, and it is also required to establish the improvability of MA estimators over least squares estimators \citep{Hansen2014risk, Zhang2016dominance}. It is noteworthy that Assumption~\ref{ass:2_add} imposes a milder constraint on $\mathcal{M}$ compared to the assumptions necessary for Theorem~\ref{theorem:1}. However, this broader applicability comes with a price: the residual of (\ref{eq:bound_class_stein}) scales with $M_n$, potentially surpassing the optimal MA risk as $M_n \to \infty$ \citep{Cavalier2001penalized}.

Given that the Stein estimators (\ref{eq:est_shr}) and (\ref{eq:stein_ma}), along with other MA estimators, are typically developed under distinct conditions on $\mathcal{M}$, it becomes crucial to set a common objective for these MA estimators for fair comparison. The subsequent subsection will address this aspect by establishing the full asymptotic optimality as a benchmark for the nested MA methods.

\subsection{Candidate model set}\label{subsec:candidate}

In this subsection, we construct the candidate model set for the Stein-type MA estimator (\ref{eq:stein_ma}) to attain the full asymptotic optimality as defined in Definition~\ref{def:fullaop}. Our approach draws inspiration from the system of weakly geometrically increasing blocks, which was initially studied in some classical nonparametric models \citep[see, e.g.,][]{Nemirovski1998, Cavalier2001penalized, Cavalier2002Sharp} and has also found application in the MA literature \citep{Dalalyan2012Sharp, Peng2023optimality}.

Specifically, consider $\mathcal{M}=\{k_1,k_2,\ldots,k_{M_n}\}$ with $k_{M_n} = p_n$. And the sizes of the candidate models in $\mathcal{M}$ are assumed to satisfy the following assumption.
\begin{assumption}\label{ass:3}
  There exists $\zeta_n > 0$ such that
  \begin{equation}\label{eq:k_ratio}
    \max_{1\leq j \leq M_n-1}\frac{k_{j+1}-k_j}{k_j-k_{j-1}} \leq 1+\zeta_n.
  \end{equation}
\end{assumption}

\begin{theorem}\label{theorem:2}
    Suppose Assumptions~\ref{ass:1}, \ref{ass:2}, and \ref{ass:3} hold. Then there exist an integer $N$ such that when $n>N$,
    \begin{equation}\label{eq:riskbound2}
  \mathbb{E}L_n(\hat{\bw} | \mathcal{M},\bmu) \leq   (1+\bar{\varphi})(1+\zeta_n) R_n(\bw^* | \mathcal{M}_a,\bmu)+\left[8c_1+k_1(1+\bar{\varphi})\right]\sigma^2,
\end{equation}
where $c_1$ is given in Assumption~\ref{ass:1}, and $\bar{\varphi}$ is defined in Theorem~\ref{theorem:1}.
\end{theorem}

This theorem indicates that the Stein-type MA estimator achieves the full asymptotic optimality under the conditions where $\bar{\varphi}\to 0$, $\zeta_n \to 0$, and the remainder term $\left[8c_1+k_1(1+\bar{\varphi})\right]\sigma^2$ is not too large. We provide several specific examples of $\mathcal{M}$ to illustrate these conditions.

Let $\nu_n$ be an integer such that $\nu_n \to \infty$ as $n \to \infty$. A typical choice is $\nu_n= \lfloor\log n \rfloor$ or $\nu_n=\lfloor \log\log n \rfloor$. And then let $\rho_n = 1/\log\nu_n$. Define a candidate model set $\mathcal{M}^*=\{k_1,\ldots,k_{M_n} \}$ with
\begin{equation*}
        k_m=\left\{\begin{array}{ll}
\nu_n &\quad m=1, \\
k_{m-1}+\lfloor  \nu_n \rho_n(1+\rho_n)^{m-1} \rfloor &\quad m=2,\ldots,M_n-1,\\
p_n &\quad m = M_n,
\end{array}\right.
      \end{equation*}
where $M_n = \arg\max_{m\in \mathbb{N}}(\nu_n+\sum_{j=2}^{m}\lfloor  \nu_n \rho_n(1+\rho_n)^{j-1} \rfloor \leq p_n)$. Then we have the following consequence.

\begin{corollary}\label{coro:1}
  Suppose
  \begin{equation}\label{eq:penalty}
    \varphi_m=\frac{1}{(k_m-k_{m-1})^{\tau}}, \quad 0<\tau<1/2,
  \end{equation}
  and
  \begin{equation}\label{eq:risk_cond}
    \nu_n = o\left[ R_n(\bw^* | \mathcal{M}_a,\bmu) \right],
  \end{equation}
  then the Stein-type MA estimator based on $\mathcal{M}^*$ is fully asymptotically optimal in terms of (\ref{eq:fullaop}).
\end{corollary}

The proof of Corollary~\ref{coro:1} involves verifying Assumptions~\ref{ass:1}, \ref{ass:2}, and \ref{ass:3} for $\mathcal{M}^*$, and it is presented in Section~\ref{sec:proof_42} of the Appendix. The hyperparameters $\tau$ and $\nu_n$ exert different influences on the Stein-type MA estimators:
\begin{itemize}
  \item For the parameter $\tau$, it does not impact the full asymptotic optimality when $0<\tau<1/2$. However, smaller $\tau$ does impede the rate at which the Stein-type MA estimator converges to the optimal MA risk. This effect stems from the slower decay of $\bar{\varphi}$ towards $0$, where $\bar{\varphi}$ contributes to the leading constant in the risk bound (\ref{eq:riskbound2}).

  \item In contrast, the parameter $\nu_n$ delineates the parameter space wherein full asymptotic optimality can be realized. For example, when $\nu_n = \lfloor\log n \rfloor$, (\ref{eq:risk_cond}) requires $R_n(\bw^* | \mathcal{M}_a,\bmu)$ to diverge faster than $\log n$ to achieve the full asymptotic optimality. Conversely, reducing the order of $\nu_n$ from $\log n$ to a slower rate expands the parameter space of the full asymptotic optimality. For example, when $\nu_n \sim \log\log n$, the Stein-type MA estimator based on $\mathcal{M}^*$ achieves the full asymptotic optimality if $R_n(\bw^* | \mathcal{M}_a,\bmu)$ diverges faster than $\log \log n$. However, it is noteworthy that employing smaller $\nu_n$ results in a larger $\rho_n$, consequently yielding larger $\zeta_n$ and thereby slowing down the convergence rate of the Stein-type MA estimator towards the optimal MA risk.
\end{itemize}

To further illustrate the above statements on the full asymptotic optimality, consider the following two representative examples with $p_n = n$.
\begin{example}[Polynomially decaying signal]

Consider $\frac{\|\bmu_{j|\mathcal{M}_a}\|}{\sqrt{n}} = j^{-\alpha_1}, 1\leq j \leq n$, where $1/2 < \alpha_1 < \infty$. Based on Theorem~1 and Example~1 of \cite{Peng2022improvability}, we know $R_n(\bw^* | \mathcal{M}_a,\bmu) \asymp n^{\frac{1}{2\alpha_1}}$. In this case, (\ref{eq:risk_cond}) is satisfied for all $1/2 < \alpha_1 < \infty$ provided $\nu_n$ increases to $\infty$ in a logarithmic rate.

\end{example}

\begin{example}[Exponentially decaying signal]

Consider $\frac{\|\bmu_{j|\mathcal{M}_a}\|}{\sqrt{n}} = \exp(-j^{\alpha_2}), 1\leq j \leq n$, where $0< \alpha_2 < \infty$. According to the results in \cite{Peng2022improvability}, we see $R_n(\bw^* | \mathcal{M}_a,\bmu) \asymp (\log n)^{\frac{1}{\alpha_2}}$. If $\nu_n \sim \log n $, (\ref{eq:risk_cond}) is satisfied when $0 < \alpha_2 < 1$. In contrast, if $\nu_n \sim \log\log n $, then (\ref{eq:risk_cond}) holds for all $0< \alpha_2 < \infty$. Therefore, the Stein-type MA estimator with $\nu_n \sim \log\log n $ achieves the full asymptotic optimality across the whole parameter space in this example.

\end{example}

\subsection{Several remarks}

\begin{remark}[Comparison with the classical Stein rule (\ref{eq:est_shr})]
  Combining Corollary~\ref{cor:classical_stein} with the relation (\ref{eq:ppc}), it is evident that if the candidate model set satisfies (\ref{eq:k_ratio}), then the risk of (\ref{eq:est_shr}) is upper bounded by
  \begin{equation}\label{eq:cla_stein_full}
  \mathbb{E}L_n(\hat{\bw}_{\text{\scriptsize STE} }|\mathcal{M},\bmu) \leq   (1+\zeta_n) R_n(\bw^* | \mathcal{M}_a,\bmu)+\left(4M_n+k_1\right)\sigma^2.
  \end{equation}
  When $\mathcal{M}^*$ is adopted, we have $M_n \lesssim \log n \log \nu_n$ \citep[Lemma 3.12 of][]{tsybakov2008introduction}, and the full asymptotic optimality of (\ref{eq:est_shr}) holds provided $\nu_n \vee( \log n\log \nu_n) = o\left[ R_n(\bw^* | \mathcal{M}_a,\bmu) \right]$, which is stronger than (\ref{eq:risk_cond}). No matter how small $\nu_n$ is, the current risk bound cannot imply the full asymptotic optimality of $\hat{\bmu}_{\hat{\bw}_{\text{\tiny STE} }|\mathcal{M}}$ over the super fast decaying signal (i.e., $\frac{\|\bmu_{j|\mathcal{M}_a}\|}{\sqrt{n}} = \exp(-j^{\alpha_2})$, $1 \leq \alpha_2 < \infty$). In contrast, the penalized blockwise Stein rule (\ref{eq:stein2}) can achieve the full asymptotic optimality in these cases.

\end{remark}

\begin{remark}[Comparison with the MMA estimator (\ref{eq:class_MMA})]
  The Stein-type MA estimators exclude the use of $\mathcal{M}_a$ due to Assumptions~\ref{ass:1}. In contrast, the MMA estimator based on $\mathcal{M}_a$ can achieve the full asymptotic optimality when the signals decay slowly. For instance, when $\frac{\|\bmu_{j|\mathcal{M}_a}\|}{\sqrt{n}} = j^{-\alpha_1}, 1/2 < \alpha_1 < \infty$, \cite{Peng2023optimality} prove that the risk of the MMA estimator is upper bounded by
  \begin{equation}\label{eq:peng}
  \begin{split}
      \mathbb{E}L_n(\hat{\bw}_{\text{\scriptsize MMA} }|\mathcal{M}_a,\bmu)& \leq R_n(\bw^* | \mathcal{M}_a,\bmu) \\
       & + C(\log n)^3 + C\left[R_n(\bw^* | \mathcal{M}_a,\bmu)\right]^{\frac{1}{2}}\left(\log n\right)^{\frac{3}{2}}.
  \end{split}
  \end{equation}
   This bound implies that $\hat{\bmu}_{\hat{\bw}_{\text{\tiny MMA} }|\mathcal{M}_a}$ achieve the full asymptotical optimality, since $R_n(\bw^* | \mathcal{M}_a,\bmu) \asymp n^{\frac{1}{2\alpha_1}}$ in this scenario. A comparison of the remainder terms in the risk bounds (\ref{eq:riskbound2}) and (\ref{eq:peng}) reveals that MMA converges to the oracle MA risk at a rate given by:
  \begin{equation}\label{eq:MMArem}
    \frac{(\log n)^{3}}{R_n(\bw^* | \mathcal{M}_a,\bmu)}\vee\left[\frac{(\log n)^{3}}{R_n(\bw^* | \mathcal{M}_a,\bmu)}\right]^{\frac{1}{2}},
  \end{equation}
  whereas the Stein-type MA estimator with $\nu_n = \lfloor\log n \rfloor$ converges in a rate of
  \begin{equation}\label{eq:steinrem}
    \frac{1}{\log\log n} \vee \frac{\log n}{R_n(\bw^* | \mathcal{M}_a,\bmu)},
  \end{equation}
  which is slower than (\ref{eq:MMArem}) when $R_n(\bw^* | \mathcal{M}_a,\bmu) \asymp n^{\frac{1}{2\alpha_1}}$. However, when the SNRs decay exponentially, MMA still relies on some reduced candidate model sets to achieve the full asymptotic optimality, as detailed in Section 4.4 of \cite{Peng2023optimality}. And it has a similar limitation to (\ref{eq:est_shr}) in failing to achieve the full asymptotic optimality when the SNRs decay extremely fast.
\end{remark}

%\begin{remark}[Regarding tightness of bounds]
%  The above comparisons are on the known upper bounds of the nested MA strategies. It is still an open question whether these upper bounds are tight.
%\end{remark}

\begin{remark}[Regarding the nested setup]
  This paper focuses on the optimal MA risk within the nested model framework. Further comparisons between nested MA strategies can be found in \cite{Fang2022}, \cite{Xu2022From}, and Section~3 of \cite{Peng2023optimality}. It is important to underscore that a fundamental assumption underlying the nested MA is the ordered importance of regressors. Without this prerequisite on regressor ordering, the application of nested MA strategies may encounter challenges, as the optimal nested MA risk $R_n(\bw^* | \mathcal{M}_a,\bmu)$ may diverge very fast. In scenarios where regressor order is unknown, there are two avenues for MA. Given that achieving optimal MA risk in the non-nested setup is generally impossible \citep[Theorem 7 of][]{Peng2023optimality}, one may consider to combine a reasonably small subset of non-nested models and construct achievable asymptotically optimal MA estimators \citep[see, e.g.,][]{Wan2010least, Zhang_2021_new}, which are useful in many practical problems. Another alternative is to give up the asymptotic optimality and resort to the minimax optimal aggregation of general candidate models, see, e.g., \cite{Yang2001Adaptive}, \cite{tsybakov2003optimal}, and \cite{Wang2014}.
\end{remark}

\section{Simulation studies}\label{sec:simu}

\subsection{Assessing the full asymptotic optimality of the Stein-type MA estimators}\label{sec:assess}

The data is simulated from the linear regression model (\ref{eq:model}), where $p_n=\lfloor 4 n^{2/3} \rfloor$, $x_{i1} = 1$, the regressor vectors $(x_{i2},\ldots,x_{ip_n})^{\top},i=1,\ldots,n$ are i.i.d. from the multivariate normal distribution with mean $\boldsymbol{0}_{p_n-1}$ and covariance $\Sigma_{(p_n-1)\times(p_n-1)}$, where $\Sigma_{ii}=1$ and $\Sigma_{ij}=0$ when $i\neq j$, and the random error terms $\epsilon_i$'s are i.i.d. from $N(0,\sigma^2)$ and are independent of ${x_{ij}}'$s. To verify the full asymptotic optimality, we consider two scenarios where the regression coefficients monotonically decay:
\begin{description}
  \item[Case 1]  $\bbeta = (\beta_1,\ldots,\beta_{p_n})^{\top}$, where $\beta_j=j^{-\alpha_1}$ and $\alpha_1$ is set to be $0.75$, $1$, and $1.5$.
  \item[Case 2] $\bbeta = (\beta_1,\ldots,\beta_{p_n})^{\top}$, where $\beta_j=\exp(-j^{\alpha_2})$ and $\alpha_2$ is set to be $0.5$, $1$, and $1.5$.
\end{description}
Let $\bbeta_{-1}$ denote $(\beta_2,\ldots,\beta_{p_n})^{\top}$. The population SNR, defined by $\bbeta_{-1}^{\top}\Sigma\bbeta_{-1}/\sigma^2$, is set to be one via the parameter $\sigma^2$. Moreover, the sample size $n$ increases from $50$ to $1500$.

Two representative Stein-type MA estimators are considered. The first (SMA1) is the penalized Stein-type MA based on the candidate model set $\mathcal{M}^*$ with $\tau=1/3$ and $\nu_n = \lfloor\log n \rfloor$. Another one (SMA2) is the same as SMA1, except that $\nu_n$ is set to be $\lfloor\log\log n \rfloor$. The candidate models used to implement MA are nested and estimated by least squares.

\begin{figure}[!t]
    \centering
        \subfigure[Case 1]{
    \begin{minipage}[t]{1\linewidth}
    \centering
       \includegraphics[width=5.5in]{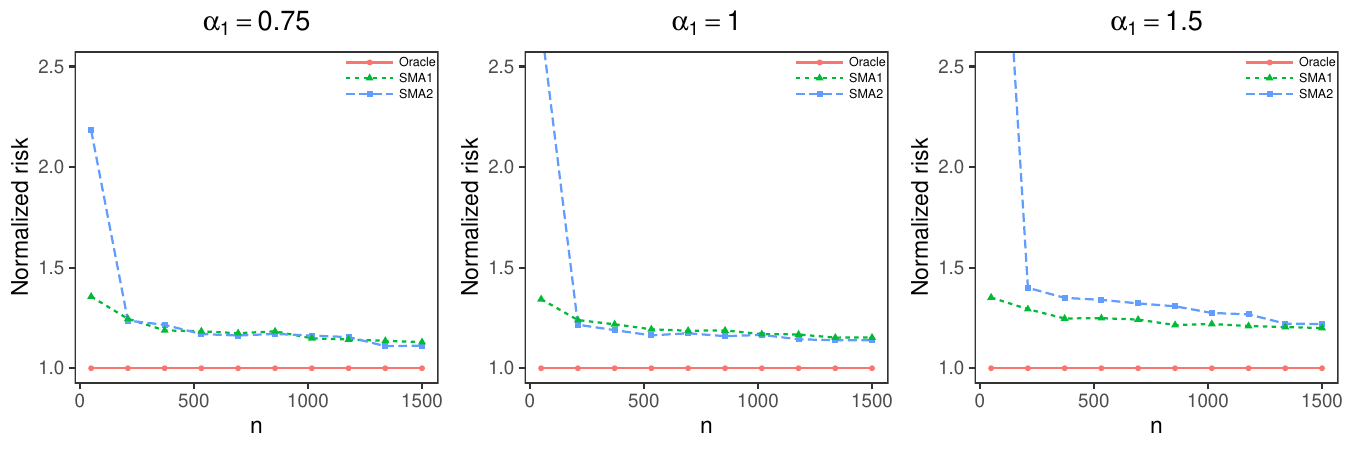}
       % \hspace{2cm}
    \end{minipage}
    }

    \subfigure[Case 2]{
    \begin{minipage}[t]{1\linewidth}
    \centering
       \includegraphics[width=5.5in]{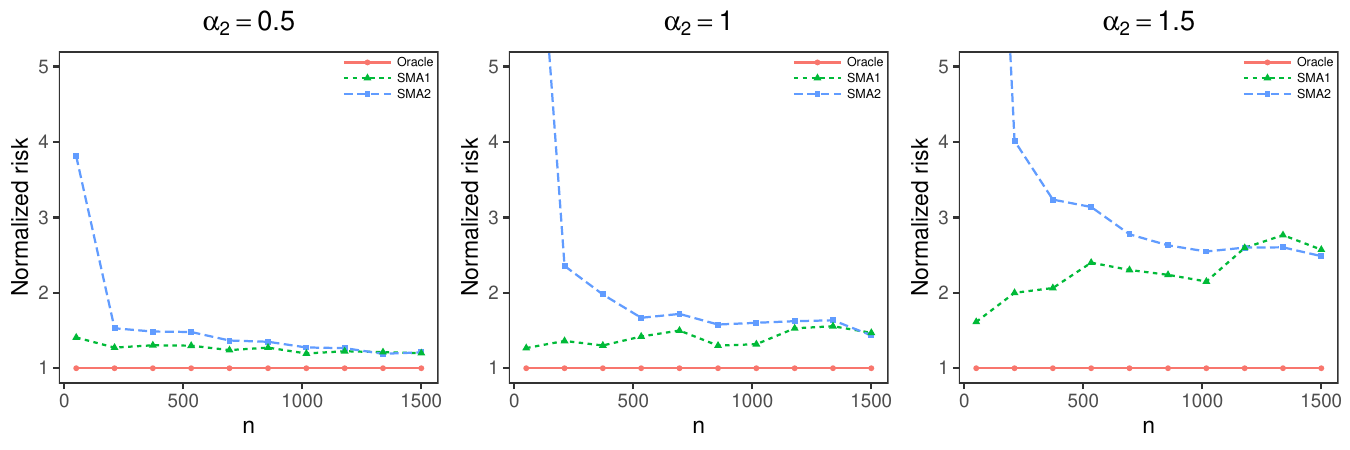}
       % \hspace{2cm}
    \end{minipage}
    }

    \caption{Normalized risks of two Stein-type MA methods with Case 1 in row (a) and Case 2 in row (b). }
    \label{fig:c1}
\end{figure}

Let $\bmu=(\mu_1,\ldots,\mu_n)^{\top}$ denote the mean vector of the true regression function. The accuracy of an estimation procedure is evaluated in terms of the squared $\ell_2$-loss $\|\bmu - \hat{\bmu} \|^2$, where $\hat{\bmu}=(\mu_1,\ldots,\mu_n)^{\top}$ is the estimated mean vector. We replicate the data generation process $R=100$ times to approximate the risks of the competing methods. Let $\bmu^{(r)}$ and $\hat{\bmu}^{(r)}$ denote the true mean vector and the estimated mean vector in the $r$-th replicate, respectively. We plot the normalized risk
\begin{equation}\label{eq:risk_ratio}
  \text{Normalized risk} = \frac{R^{-1}\sum_{r=1}^{R}\|\bmu^{(r)} - \hat{\bmu}^{(r)} \|^2}{\min_{\mathbf{w}\in\mathcal{W}_{p_n}}\left\{R^{-1}\sum_{r=1}^{R}\|\bmu^{(r)} - \hat{\bmu}^{(r)}_{\mathbf{w}|\mathcal{M}_a} \|^2\right\}}
\end{equation}
as a function of $n$, where the denominator of (\ref{eq:risk_ratio}) is an estimate of the oracle risk of MA based on $\mathcal{M}_a$. The simulation results are displayed in Figure~\ref{fig:c1}.

From Figure~\ref{fig:c1} (a), we observe that the normalized risks of the methods SMA1 and SMA2 gradually decrease and approach 1 with increasing $n$, which supports the full asymptotic optimality statements in Section~\ref{subsec:candidate}. The patterns presented in Figure~\ref{fig:c1} (b) are also consistent with the established theories. When $\alpha_2=0.5$, there are slight decreases in the normalized risks of the Stein-type estimator SMA1. In contrast, the curves of SMA1 level off with relative risks significantly larger than 1 when $\alpha_2 \geq 1$. This result corroborates the upper bounds in Section~\ref{subsec:candidate} that SMA1 achieves the full asymptotic optimality when $\alpha_2<1$ but may not when the coefficients decay too fast. When $\alpha_2=1.5$, it can be seen from Figure~\ref{fig:c1} (b) that the normalized risk of SMA2 declines steadily from 14.46 to 2.48, whereas the curve of SMA1 increases slightly from 1.61 to 2.57. This finding is consistent with our theoretical understanding in Corollary~\ref{coro:1} that SMA2 can achieve the fully asymptotic optimality in a broader parameter space than SMA1. On the other hand, SMA2 actually performs worse than SMA1, especially when $n\leq 500$. Thus, it also illustrates that the asymptotic optimality should not be the sole justification for a MA method. As shown in Figure~\ref{fig:c1} (b), the estimator with satisfactory asymptotic property (SMA2) may converge slowly and perform poorly in the finite sample setting.

\subsection{Comparing different MA estimators}

In this subsection, we investigate the finite sample performance of several MA strategies based on different candidate model sets and weight sets. Besides SMA1 and SMA2 discussed in Section~\ref{sec:assess}, we consider additional five competing methods for comparison. The first (MMA1) is the classical MMA method, which minimizes (\ref{eq:criterion}) over $\mathcal{W}_{p_n}$ to combine all nested models. The second (MMA2) is a parsimonious version of the MMA method \citep{Zhang2020Parsimonious}, which uses the same weight set and candidate model set as MMA1 but replaces $2\sigma^2$ in (\ref{eq:criterion}) with a stronger penalty factor $\log n$. The third MMA method (MMA3) employs the same candidate model set as SMA1 but minimizes the MMA criterion over the standard weight constraint (\ref{eq:weight_general}). The fourth method, MMA4, utilizes the standard weight constraint and a candidate set comprising equal-sized blocks with $4$ regressors in each block. Additionally, we include a James-Stein estimator (SMA3), as investigated in Section 6 of \cite{Hansen2014risk}, which combines only the intercept model and the full model.

The simulation data are generated from the same regression model in Section~\ref{sec:assess}, with SNRs set to two. Additionally, we introduce two cases of monotonically increasing coefficients, which may be less favorable for the nested MA with the standard weight constraints.
\begin{description}
  \item[Case 3]   $\bbeta = (\beta_1,\ldots,\beta_{p_n})^{\top}$, where $\beta_j=(p_n+1-j)^{-\alpha_3}$ and $\alpha_3$ is set to be $0.75$, $1$, and $1.5$.
  \item[Case 4]   $\bbeta = (\beta_1,\ldots,\beta_{p_n})^{\top}$, where $\beta_j=\exp[-(p_n+1-j)^{\alpha_4}]$ and $\alpha_4$ is set to be $0.5$, $1$, and $1.5$.
\end{description}
To better assess the performance of the competing methods and uncover additional underlying properties beyond those depicted in Figure~\ref{fig:c1}, we divide the risk of each method by the optimal MA risk based on the candidate model set $\mathcal{M}^*$ with $\nu_n = \lfloor\log n \rfloor$ and the weight set $\mathcal{W}_{|\mathcal{M}^*|}$. The simulation results are presented in Tables~\ref{tab:t1}--\ref{tab:t2} (the standard errors are in the parentheses).

% Table generated by Excel2LaTeX from sheet 'Sheet2'
\begin{table}[!t]
  \centering
  \caption{Comparisons of different MA estimators in Cases~1--2. The risk of each method is divided by the optimal MA risk based on $\mathcal{M}^*$ and $\mathcal{W}_{|\mathcal{M}^*|}$ in each simulation. }
  \resizebox{\columnwidth}{!}{
    \begin{tabular}{llcccccc}
    \hline
       \multirow{2}[0]{*}{$n$}   &   \multirow{2}[0]{*}{method}     & \multicolumn{3}{c}{Case 1} & \multicolumn{3}{c}{Case 2} \\
       \cline{3-5} \cline{6-8}
      & & \multicolumn{1}{c}{$\alpha_1$=0.75} & \multicolumn{1}{c}{$\alpha_1$=1} & \multicolumn{1}{c}{$\alpha_1$=1.5} & \multicolumn{1}{c}{$\alpha_2$=0.5} & \multicolumn{1}{c}{$\alpha_2$=1} & \multicolumn{1}{c}{$\alpha_2$=1.5} \\
      \hline
    100   & MMA1  & 1.153 (0.015) & 1.244 (0.017) & 1.485 (0.044) & 1.434 (0.049) & 2.076 (0.132) & 1.961 (0.204) \\
          & MMA2  & 1.541 (0.026) & 1.552 (0.028) & 1.542 (0.040) & 1.645 (0.043) & 1.852 (0.101) & 1.147 (0.113) \\
          & MMA3  & 1.162 (0.013) & 1.191 (0.014) & 1.372 (0.038) & 1.364 (0.048) & 1.485 (0.075) & 1.793 (0.165) \\
          & MMA4  & 1.142 (0.015) & 1.212 (0.017) & 1.442 (0.043) & 1.384 (0.049) & 1.689 (0.104) & 2.059 (0.194) \\
          & SMA1  & 1.281 (0.016) & 1.306 (0.021) & 1.341 (0.030) & 1.387 (0.038) & 1.160 (0.031) & 1.171 (0.039) \\
          & SMA2  & 1.374 (0.022) & 1.368 (0.023) & 1.480 (0.042) & 1.339 (0.040) & 4.376 (0.309) & 4.992 (0.586) \\
          & SMA3  & 1.947 (0.035) & 2.883 (0.061) & 6.196 (0.227) & 5.594 (0.218) & 23.799 (1.625) & 27.888 (3.902) \\[3pt]
    500   & MMA1  & 1.094 (0.004) & 1.163 (0.008) & 1.360 (0.024) & 1.364 (0.027) & 1.786 (0.111) & 1.459 (0.146) \\
          & MMA2  & 1.825 (0.020) & 1.732 (0.021) & 1.790 (0.041) & 1.782 (0.042) & 1.411 (0.063) & 0.884 (0.097) \\
          & MMA3  & 1.080 (0.004) & 1.134 (0.007) & 1.266 (0.022) & 1.201 (0.012) & 1.730 (0.101) & 1.888 (0.234) \\
          & MMA4  & 1.086 (0.004) & 1.146 (0.008) & 1.321 (0.023) & 1.232 (0.016) & 2.254 (0.130) & 1.768 (0.242) \\
          & SMA1  & 1.172 (0.009) & 1.211 (0.011) & 1.288 (0.022) & 1.271 (0.021) & 1.134 (0.024) & 1.218 (0.044) \\
          & SMA2  & 1.156 (0.007) & 1.193 (0.009) & 1.416 (0.036) & 1.427 (0.027) & 1.743 (0.091) & 2.012 (0.209) \\
          & SMA3  & 2.142 (0.016) & 3.973 (0.049) & 11.337 (0.275) & 11.397 (0.324) & 48.596 (3.083) & 56.233 (6.425) \\[3pt]
    1000  & MMA1  & 1.065 (0.003) & 1.143 (0.006) & 1.311 (0.019) & 1.277 (0.016) & 2.076 (0.118) & 1.699 (0.237) \\
          & MMA2  & 1.893 (0.016) & 1.804 (0.021) & 1.788 (0.033) & 1.811 (0.034) & 1.722 (0.096) & 0.789 (0.055) \\
          & MMA3  & 1.057 (0.003) & 1.110 (0.005) & 1.240 (0.016) & 1.004 (0.001) & 1.003 (0.001) & 1.004 (0.001) \\
          & MMA4  & 1.061 (0.003) & 1.130 (0.006) & 1.277 (0.019) & 1.003 (0.001) & 1.002 (0.001) & 1.002 (0.001) \\
          & SMA1  & 1.135 (0.005) & 1.176 (0.009) & 1.253 (0.015) & 1.224 (0.014) & 1.241 (0.046) & 1.254 (0.082) \\
          & SMA2  & 1.174 (0.006) & 1.181 (0.008) & 1.233 (0.015) & 1.214 (0.015) & 1.514 (0.057) & 1.540 (0.078) \\
          & SMA3  & 2.138 (0.014) & 4.530 (0.051) & 13.841 (0.281) & 15.832 (0.402) & 83.044 (5.459) & 85.257 (6.571) \\
          \hline
    \end{tabular}%
    }
  \label{tab:t1}%
\end{table}%

From Table~\ref{tab:t1}, we see the relative risks of MMA3 and SMA1 gradually decease and approach 1 in Case 1, which support the asymptotic optimality of these two methods in the polynomially decaying cases. However, in situations where the coefficients decay fast, as in Case 2, the relative risks of MMA3 continue to exhibit an obvious downward trend. In contrast, SMA1 demonstrates a slower convergence compared to MMA3 (specifically when $\alpha_2 = 0.5$) and even displays a slight increase in risk when the coefficients decay extremely fast (i.e., $\alpha_2 = 1.5$). Although Theorem~\ref{theorem:1} asserts that SMA1 can asymptotically achieve $R_n(\bw^* | \mathcal{M}^*,\bmu)$, it is noteworthy that, in comparison to MMA3, SMA1 minimizes the principle of URE within a larger weight set. As demonstrated in Theorem~\ref{theo:improv}, this weight relaxation does not provide benefit when the coefficients are monotonically decreasing. Therefore, SMA1 may incur a higher cost to attain the optimal MA risk compared to MMA3, and converges slower especially when the coefficients decay fast. According to the similar reason, the MMA1 and MMA4 methods, although combine more candidate models than SMA1 and SMA2, may also show advantage in the finite sample setting in Cases~1--2, respectively.

As observed in Table~\ref{tab:t1}, the relative risks of the MMA2 method do not decrease in Case 1 and have values around 1.8. This result may be explained by the fact that the MMA2 method tends to assign more weights to small models due to its strong penalty, while the optimal size of candidate model in Case 1 is relatively large. When $\alpha_2=1.5$, the coefficients decay extremely fast, and MA does not have any real benefit compared to MS \citep{Peng2022improvability}. Such a case is close to the setting where a fixed-dimension true model exists. In this case, the MMA2 method becomes the best when $n=1000$. This result is expected and in accord with the theory in \cite{Zhang2020Parsimonious}. The SMA3 method, which exclusively combines the intercept model and full model, exhibit subpar performance across both cases due to a poor trade-off between bias and variance.

% Table generated by Excel2LaTeX from sheet 'Sheet2'
\begin{table}[!t]
  \centering
  \caption{Comparisons of different MA estimators in Cases~3--4. The risk of each method is divided by the optimal MA risk based on $\mathcal{M}^*$ and $\mathcal{W}_{|\mathcal{M}^*|}$ in each simulation. }
  \resizebox{\columnwidth}{!}{
    \begin{tabular}{llcccccc}
    \hline
         \multirow{2}[0]{*}{$n$}   &   \multirow{2}[0]{*}{method}      & \multicolumn{3}{c}{Case 3} & \multicolumn{3}{c}{Case 4} \\
         \cline{3-5} \cline{6-8}
        &  & \multicolumn{1}{c}{$\alpha_3=0.75$} & \multicolumn{1}{c}{$\alpha_3=1$} & \multicolumn{1}{c}{$\alpha_3=1.5$} & \multicolumn{1}{c}{$\alpha_4=0.5$} & \multicolumn{1}{c}{$\alpha_4=1$} & \multicolumn{1}{c}{$\alpha_4=1.5$} \\
        \hline
    100   & MMA1  & 1.029 (0.004) & 1.035 (0.004) & 1.050 (0.006) & 1.044 (0.005) & 1.034 (0.005) & 1.028 (0.004) \\
          & MMA2  & 1.806 (0.034) & 1.797 (0.034) & 1.810 (0.043) & 1.881 (0.041) & 1.836 (0.039) & 1.861 (0.034) \\
          & MMA3  & 1.035 (0.003) & 1.034 (0.003) & 1.050 (0.005) & 1.046 (0.005) & 1.040 (0.005) & 1.033 (0.003) \\
          & MMA4  & 1.043 (0.003) & 1.047 (0.004) & 1.062 (0.006) & 1.056 (0.005) & 1.046 (0.005) & 1.040 (0.004) \\
          & SMA1  & 1.111 (0.011) & 1.107 (0.011) & 1.127 (0.013) & 1.137 (0.012) & 1.121 (0.013) & 1.126 (0.010) \\
          & SMA2  & 1.055 (0.009) & 1.054 (0.009) & 1.076 (0.011) & 1.080 (0.010) & 1.073 (0.012) & 1.064 (0.007) \\
          & SMA3  & 1.007 (0.004) & 1.008 (0.004) & 1.022 (0.006) & 1.018 (0.005) & 1.014 (0.005) & 1.01 (0.004) \\[3pt]
    500   & MMA1  & 0.999 (0.001) & 1.002 (0.001) & 1.001 (0.001) & 1.002 (0.001) & 1.001 (0.001) & 1.000 (0.001) \\
          & MMA2  & 2.104 (0.023) & 2.180 (0.028) & 2.132 (0.025) & 2.120 (0.025) & 2.136 (0.024) & 2.155 (0.024) \\
          & MMA3  & 1.006 (0.001) & 1.008 (0.001) & 1.007 (0.001) & 1.008 (0.001) & 1.007 (0.001) & 1.007 (0.001) \\
          & MMA4  & 1.004 (0.001) & 1.006 (0.001) & 1.004 (0.001) & 1.006 (0.001) & 1.005 (0.001) & 1.004 (0.001) \\
          & SMA1  & 1.015 (0.005) & 1.027 (0.007) & 1.021 (0.005) & 1.016 (0.006) & 1.011 (0.005) & 1.019 (0.006) \\
          & SMA2  & 1.003 (0.004) & 1.013 (0.005) & 1.012 (0.005) & 1.000 (0.004) & 1.002 (0.004) & 1.009 (0.004) \\
          & SMA3  & 0.998 (0.001) & 1.000 (0.001) & 1.000 (0.001) & 1.000 (0.001) & 1.000 (0.001) & 0.998 (0.001) \\[3pt]
    1000  & MMA1  & 1.001 (0.001) & 1.001 (0.001) & 1.000 (0.001) & 1.001 (0.001) & 1.000 (0.001) & 1.001 (0.001) \\
          & MMA2  & 2.204 (0.020) & 2.192 (0.018) & 2.211 (0.017) & 2.211 (0.021) & 2.193 (0.018) & 2.234 (0.021) \\
          & MMA3  & 1.004 (0.001) & 1.004 (0.001) & 1.003 (0.001) & 1.004 (0.001) & 1.004 (0.001) & 1.004 (0.001) \\
          & MMA4  & 1.002 (0.001) & 1.003 (0.001) & 1.002 (0.001) & 1.003 (0.001) & 1.002 (0.001) & 1.003 (0.001) \\
          & SMA1  & 1.000 (0.005) & 0.994 (0.004) & 0.985 (0.004) & 0.985 (0.004) & 0.990 (0.004) & 0.993 (0.005) \\
          & SMA2  & 1.005 (0.003) & 1.000 (0.003) & 0.996 (0.003) & 0.998 (0.003) & 1.001 (0.002) & 1.003 (0.003) \\
          & SMA3  & 1.000 (0.001) & 1.000 (0.001) & 0.999 (0.001) & 1.000 (0.001) & 1.000 (0.001) & 1.000 (0.001) \\
          \hline
    \end{tabular}%
    }
  \label{tab:t2}%
\end{table}%

The results presented in Table~\ref{tab:t2} are different from those in Table~\ref{tab:t1}, yet remain consistent with the current MA theories. Specifically, when the coefficients monotonically increase, the optimal weight vector within (\ref{eq:weight_general}) typically allocates weight to the intercept and full models. Consequently, the relative risks of MMA1, MMA3, and MMA4—minimizing the MMA criterion over the standard weight constraints—converge to 1 as $n$ increases. The relative risks of the James-Stein estimator SMA3, which only combines the intercept model and the full model, rapidly approach 1. However, in such scenarios, the optimal MA risk within the standard constraint does not represent the optimal performance. The relative risks of SMA1 and SMA2 fall below one across various cases when $n=1000$, owing to their targets on the optimal risks within the relaxed weight sets. As suggested by the results in Section~\ref{subsec:relax}, the weight relaxation may have advantage when the coefficients are not monotonically non-increasing.

Moreover, it is illuminating to observe that the advantage of SMA1 and SMA2 is not substantial, and in certain cases (e.g., $\alpha_4=1.5$), they may not have tangible improvement. This limitation arises from the construction of $\mathcal{M}^*$, which is based on a sequence of increasing-size blocks. This construction actually assumes that the regressors are ordered from the most important to least important. The above analysis also implies that while SMA1 and SMA2 are based on the relaxed weight constraints, they inherently favor monotonic non-increasing signals and therefore are inadequate for the scenarios where the signal's monotonicity is ambiguous or unknown.

\subsection{Comparing different choices of $\tau$ and $\nu_n$}

\begin{figure}[!t]
    \centering
        \subfigure[$n=500$]{
    \begin{minipage}[t]{1\linewidth}
    \centering
       \includegraphics[width=4.7in]{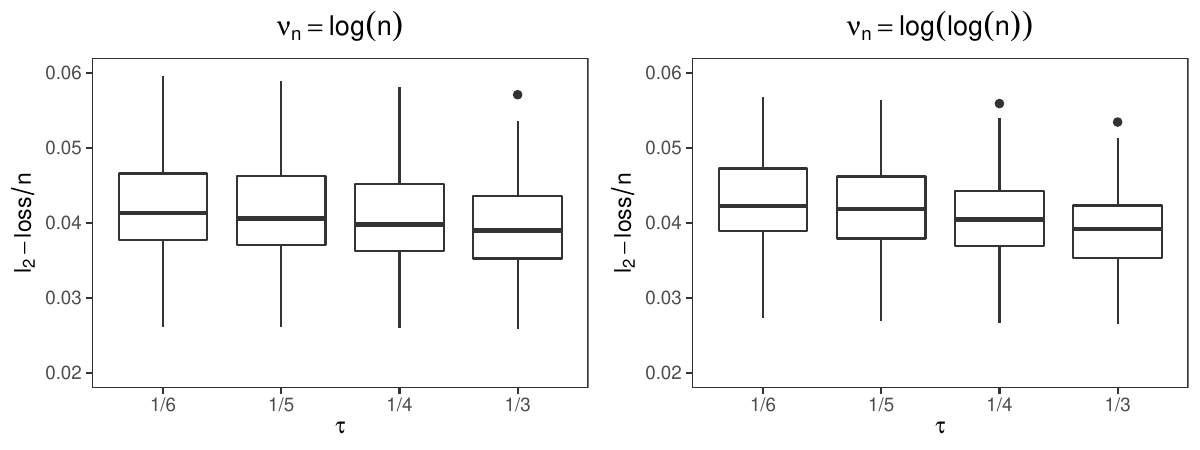}
       % \hspace{2cm}
    \end{minipage}
    }

    \subfigure[$n=1000$]{
    \begin{minipage}[t]{1\linewidth}
    \centering
       \includegraphics[width=4.7in]{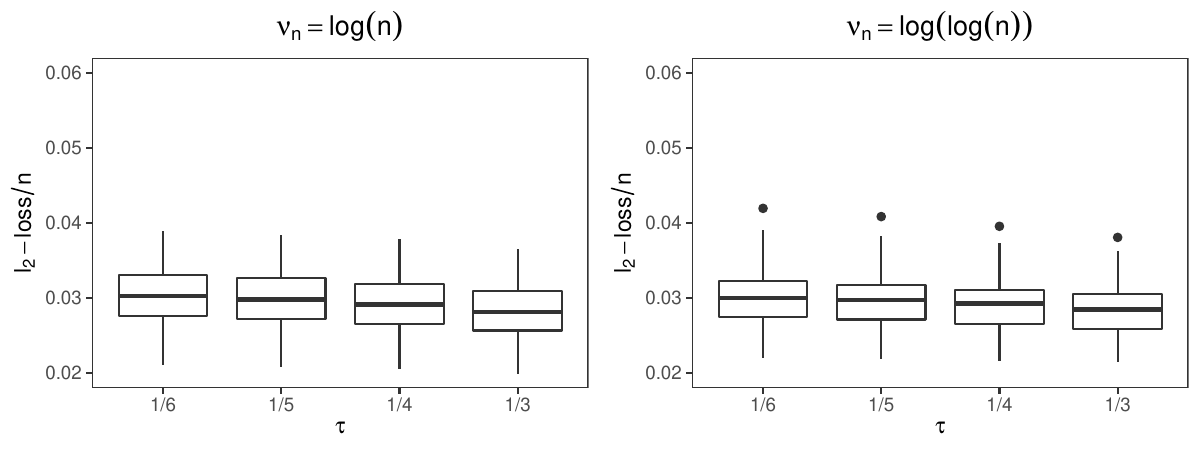}
       % \hspace{2cm}
    \end{minipage}
    }

    \caption{Comparisons of different choices of $\tau$ and $\nu_n$ when $\beta_j=j^{-1}$.}
    \label{fig:c2}
\end{figure}

\begin{figure}[!t]
    \centering
        \subfigure[$n=500$]{
    \begin{minipage}[t]{1\linewidth}
    \centering
       \includegraphics[width=4.7in]{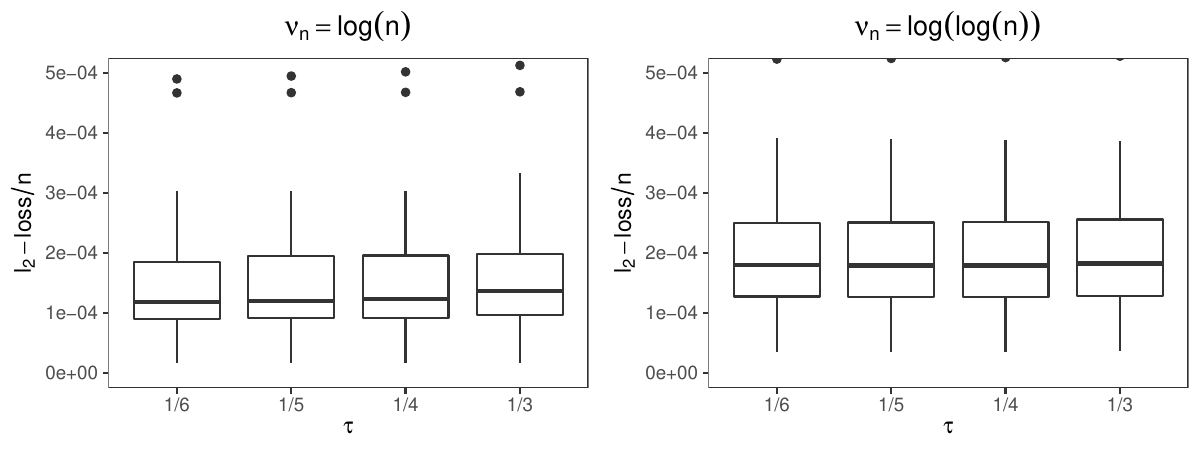}
       % \hspace{2cm}
    \end{minipage}
    }

    \subfigure[$n=1000$]{
    \begin{minipage}[t]{1\linewidth}
    \centering
       \includegraphics[width=4.7in]{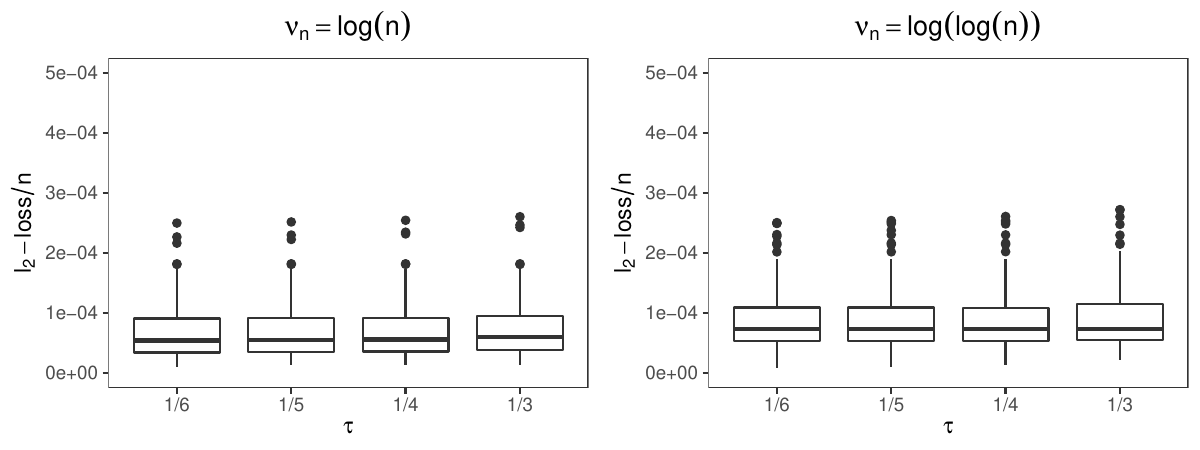}
       % \hspace{2cm}
    \end{minipage}
    }

    \caption{Comparisons of different choices of $\tau$ and $\nu_n$ when $\beta_j=\exp(-j)$.}
    \label{fig:c3}
\end{figure}

In this subsection, we investigate several different choices of the hyperparameters $\tau$ and $\nu_n$ in the Stein-type MA estimator (\ref{eq:stein_ma}) based on $\mathcal{M}^*$. The simulation results are displayed in Figures~\ref{fig:c2}--\ref{fig:c3}.

Overall, it is observed that the performance of the Stein-type MA estimator is robust across different parameter choices when the coefficients decay polynomially. The variations of the mean values in each panel of Figure~\ref{fig:c2} are not significant. For example, when $\beta_j = j^{-1}$, $n=500$, and $\nu_n=\log n$, the p-value from a t-test of comparing $\tau = 1/6$ with $\tau=1/3$ is 0.0058, indicating no significant difference in this case. When $\tau=1/3$, the p-value of comparing $\nu_n=\log n$ v.s. $\nu_n=\log\log n$ is 0.0011. However, it is noted that the performance of the Stein-type estimator becomes more sensitive to the hyperparameters when the coefficients decay exponentially. For example, when $\beta_j = \exp(-j)$, $n=500$, and $\nu_n=\log n$, the p-value of comparing $\tau = 1/6$ with $\tau=1/3$ is 0.3654, indicating a more noticeable impact of $\tau$ on performance in this case. %The above phenomenon is expected. Since $\theta_j = j^{-\alpha_1}$, the optimal nested models scale with $n^{1/(2\alpha_1)}$, the changes of the hyperparameter only changes the sizes of candiate models in a logarithm factor, hence does not impact too much. In contrast, in the case of $\theta_j = \exp(-j^{\alpha_1})$, the optimal sizes also scale with a logarithm scale. In this case, the changes of hyperparameters may be sensitive.

Although MA is initially proposed to mitigate MS uncertainty, the above observations suggest that it is still necessary to apply some data-driven techniques to combine or select MA methods for better adaptivity, especially in real world applications. For further exploration on this topic, see \cite{ZHANG201595, QIAN2022193, FANG2022110683}.

\section{Discussion}\label{sec:discs}

This paper establishes an explicit link between MA and shrinkage in a multiple model setting, which significantly enhances the previous understanding of the relationship between MA and shrinkage in the two-model settings. It is revealed that the MMA estimator can be viewed as a variant of the positive-part Stein estimator, as both are derived from the principle of URE. The key distinction lies in the optimization approach: MMA minimizes the principle of URE within a unit simplex, whereas the latter operates within a more relaxed weight set. Building upon the established connections, we extend the penalized blockwise Stein rule to the linear regression setting to develop the asymptotically optimal MA estimators. We provide some specific candidate model sets on which the proposed Stein-type MA estimator achieves the performance of the optimal convex combination of all the nested models (i.e., the full asymptotic optimality). The improvement of the proposed Stein-type MA over the existing MA approaches is illustrated theoretically and numerically. Note that a limitation of our Stein-type MA method is that it requires the variance of the error terms to be known. Thus, extending our results to the case of unknown variance is a pressing topic for future research.

%A remarkable observation from our simulation results is that the asymptotical optimality can never be the sole justification for an MA procedure. Indeed, an asymptotically optimal method with the widest parameter region (such as SMA2 in our simulation) performs much worse when the sample size is small. In practice, it is more important to compare different MA methods according to their finite sample properties.

The unveiled connections between MA and shrinkage offer the possibility of novel methodological developments in the area of MA. The focus of this paper has been on a linear regression setting. It is of great interest to bridge the gap between MA and shrinkage in the generalized linear model setting, and then apply the Stein estimators in some general distribution families \citep[e.g., see Chapter 5 of][for a review]{Hoffmann2000} to combine models. In addition, given the approximate and exact distributions of the Stein-type estimators \citep{ULLAH1982305, Phillips1984} and the techniques of constructing Stein-type confidence intervals \citep{Hwang1982, He1992}, it is greatly desirable to conduct inference for the asymptotically optimal MA estimators. Note that \cite{Hansen2014risk} and \cite{Zhang_liu_2019} have previously investigated the inference of MA but without the asymptotic optimality results. Another research direction is building a unified theory for Bayesian and frequentist MA. Indeed, the BIC weighting method considered in the frequentist literature \citep{Buckland1997, Hjort2003} can be seen as an approximation of Bayesian MA. We conjecture that the asymptotically optimal MA estimator may also have a Bayesian interpretation since the Stein-type estimation is essentially an empirical Bayes approach \citep[see, e.g.,][]{Efron1973Stein}. We leave these for future work.

%%%%%%%%%%%%%%%%%%%%%%%%%%%%%%%%%%%%%%%%%%%%%%
%% Example with single Appendix:            %%
%%%%%%%%%%%%%%%%%%%%%%%%%%%%%%%%%%%%%%%%%%%%%%

\begin{appendix}

\renewcommand{\thesection}{A.\arabic{section}}
\numberwithin{equation}{section}

\section{Proof of Theorem~\ref{theo:improv}}\label{sec:proof:improv}

From (\ref{eq:risk}), we have
\begin{equation}\label{eq:addd3}
  \begin{split}
R_n(\bw | \mathcal{M}, \bmu)& =R_n(\bgamma | \mathcal{M}, \bmu)\\
  & = \sum_{m=1}^{M_n}\left(\|\bmu_{m|\mathcal{M}}\|^2+\sigma^2_{m|\mathcal{M}}\right)\left(\gamma_m - \frac{\|\bmu_{m|\mathcal{M}}\|^2}{\|\bmu_{m|\mathcal{M}}\|^2+\sigma^2_{m|\mathcal{M}}}  \right)^2\\
   &+ \sum_{m=1}^{M_n}\frac{\|\bmu_{m|\mathcal{M}}\|^2\sigma^2_{m|\mathcal{M}}}{\|\bmu_{m|\mathcal{M}}\|^2+\sigma^2_{m|\mathcal{M}}}+ \|\bmu_{M_n+1|\mathcal{M}}\|^2.
  \end{split}
\end{equation}
Since $0\leq\|\bmu_{m|\mathcal{M}}\|^2/(\|\bmu_{m|\mathcal{M}}\|^2+\sigma^2_{m|\mathcal{M}})\leq1$, the optimal MA risk in the enlarged weight set $\tilde{\mathcal{W}}_{M_n}$ is
\begin{equation*}
\begin{split}
   R_n(\tilde{\bw}^*|\mathcal{M}, \bmu) &= \min_{\bw \in \tilde{\mathcal{W}}_{M_n}}R_n(\bw | \mathcal{M}, \bmu)=\min_{\bgamma \in \tilde{\Gamma}_{M_n}}R_n(\bgamma | \mathcal{M}, \bmu) \\
     & = \sum_{m=1}^{M_n}\frac{\|\bmu_{m|\mathcal{M}}\|^2\sigma^2_{m|\mathcal{M}}}{\|\bmu_{m|\mathcal{M}}\|^2+\sigma^2_{m|\mathcal{M}}}+\|\bmu_{M_n+1|\mathcal{M}}\|^2,
\end{split}
\end{equation*}
which proves (\ref{eq:optimal_enlarge}).

The proof of (\ref{eq:optimal_convex}) involves a list of notations which are defined in an iterative manner. In order to ease the reader into the proof, we summarize these notations in  Algorithm~\ref{alg:ALG1}. We also outline in plain English how the first few iterations work.
\begin{algorithm}
        \caption{The notations used to prove (\ref{eq:optimal_convex})}
        \label{alg:ALG1}
        \begin{algorithmic}
        \Require $\mathcal{M}=\{k_1,k_2,\ldots,k_{M_n} \}$, $\Gamma = \left\{\bgamma :1=\gamma_1\geq \gamma_2 \geq \cdots \geq \gamma_{M_n}\geq 0\right\}$
        \Ensure $\mathcal{M}_{T}$
        \State $t \gets 1$
        \State $\mathcal{M}_{t} \gets \mathcal{M}$
        \State $\Gamma_t \gets \Gamma$
        \State $\mathcal{I}_t \gets \emptyset$
        \State $\mathcal{J}_t \gets \emptyset$
        \While{$|\mathcal{M}_t|>2$ \text{\textbf{and}} $\exists\, l \in \{2,\ldots,|\mathcal{M}_t|-1 \}$ \text{ s.t. } $\frac{\|\bmu_{l|\mathcal{M}_t}\|^2}{\sigma^2_{l|\mathcal{M}_t}}< \frac{\|\bmu_{l+1|\mathcal{M}_t}\|^2}{\sigma^2_{l+1|\mathcal{M}_t}}$}
        \State $\mathcal{I}_{t+1} \gets \left\{ l \in \left\{2,\ldots,|\mathcal{M}_t|-1 \right\}:\frac{\|\bmu_{l|\mathcal{M}_t}\|^2}{\sigma^2_{l|\mathcal{M}_t}}< \frac{\|\bmu_{l+1|\mathcal{M}_t}\|^2}{\sigma^2_{l+1|\mathcal{M}_t}}  \right\}$
        \State $\mathcal{M}_{t+1} \gets \mathcal{M}_t\setminus \{\mathcal{M}_t(l):l\in  \mathcal{I}_{t+1}\}$ \Comment{$\mathcal{M}_t(l)$ is the $l$-th smallest element in $\mathcal{M}_t$.}
        \State $\mathcal{J}_{t+1} \gets \mathcal{J}_{t}\cup\{m: k_m \in \{\mathcal{M}_t(l):l\in  \mathcal{I}_{t+1}\}  \}$
        \State $\Gamma_{t+1} \gets \left\{\bgamma\in \Gamma : \gamma_l = \gamma_{l+1},l\in \mathcal{J}_{t+1}\right\}$
        \State $t \gets t +1$
        \EndWhile
        \State $ T \gets t $
        \end{algorithmic}
\end{algorithm}

In this part, we consider the nested model set $\mathcal{M}=\{k_1,k_2,\ldots,k_{M_n} \}$, where $1\leq k_1 \leq k_2 \leq \cdots \leq k_{M_n} \leq p_n$. For any candidate model set $\mathcal{A}$, we define $\mathcal{A}(l)$ as the $l$-smallest model size in $\mathcal{A}$. For example, $\mathcal{M}(l)=k_l$ for $1\leq l \leq M_n$. In addition,
we write $\Gamma = \Gamma_{M_n}$ for the simplicity of notation and assume $M_n\geq 3$.

If $\|\bmu_{l|\mathcal{M}}\|^2/\sigma^2_{l|\mathcal{M}}\geq \|\bmu_{l+1|\mathcal{M}}\|^2/\sigma^2_{l+1|\mathcal{M}}$ for $2 \leq l \leq M_n-1$, then we have $T=1$ and $\mathcal{M}_T=\mathcal{M}_1=\mathcal{M}$. In view of (\ref{eq:addd3}), the optimal MA risk in this case is given by
\begin{equation*}
\begin{split}
R_n(\bw^*|\mathcal{M}, \bmu) &= \min_{\bw \in \mathcal{W}_{M_n}}R_n(\bw | \mathcal{M}, \bmu) = \min_{\bgamma \in \Gamma}R_n(\bgamma | \mathcal{M}, \bmu) \\
     & = \frac{\sigma^4_{1|\mathcal{M}}}{\|\bmu_{1|\mathcal{M}}\|^2+\sigma^2_{1|\mathcal{M}}}+\sum_{m=2}^{M_n}\frac{\|\bmu_{m|\mathcal{M}}\|^2\sigma^2_{m|\mathcal{M}}}{\|\bmu_{m|\mathcal{M}}\|^2+\sigma^2_{m|\mathcal{M}}}+\|\bmu_{M_n+1|\mathcal{M}}\|^2\\
     & = R_n(\bw^*|\mathcal{M}_T, \bmu),
\end{split}
\end{equation*}
where the third equality follows from that $\gamma_1=1$ and $\|\bmu_{m|\mathcal{M}}\|^2/(\|\bmu_{m|\mathcal{M}}\|^2+\sigma^2_{m|\mathcal{M}})$ are monotonically non-increasing, and the last equality is due to $\mathcal{M} = \mathcal{M}_T$. Thus, (\ref{eq:optimal_convex}) is proved.

When the sequence $\{ \|\bmu_{m|\mathcal{M}}\|^2/\sigma^2_{m|\mathcal{M}}\}_{m=2}^{M_n}$ is not strictly monotonically non-increasing, let us process the first iteration in Algorithm~\ref{alg:ALG1} with $t=1$ and $\mathcal{M}_1=\mathcal{M}$. Let
\begin{equation*}
  \begin{split}
     \mathcal{I}_2 & \triangleq \left\{ l \in \left\{2,\ldots,|\mathcal{M}_1|-1 \right\}:\|\bmu_{l|\mathcal{M}_1}\|^2/\sigma^2_{l|\mathcal{M}_1}< \|\bmu_{l+1|\mathcal{M}_1}\|^2/\sigma^2_{l+1|\mathcal{M}_1}  \right\} \\
       & = \left\{ l \in \left\{2,\ldots,M_n-1 \right\}:\|\bmu_{l|\mathcal{M}}\|^2/\sigma^2_{l|\mathcal{M}}< \|\bmu_{l+1|\mathcal{M}}\|^2/\sigma^2_{l+1|\mathcal{M}}  \right\}
  \end{split}
\end{equation*}
denote the set of model index in $\mathcal{M}_1$ where the monotonicity is violated. Define
\begin{equation*}
\begin{split}
   \mathcal{J}_2 & \triangleq \left\{m: k_m \in \{\mathcal{M}_1(l):l\in  \mathcal{I}_2\}  \right\}  \\
     & = \left\{m: k_m \in \{k_l: l\in  \mathcal{I}_2\}  \right\}\\
     & = \mathcal{I}_2,
\end{split}
\end{equation*}
where the second equality follows from $\mathcal{M}_1= \mathcal{M}$. In addition, define
\begin{equation*}\label{eq:monotone_1}
  \Gamma_{2}\triangleq\left\{\bgamma\in \Gamma : \gamma_l = \gamma_{l+1},\,l\in \mathcal{J}_2\right\}
\end{equation*}
as the subset of $\Gamma$ with the equality constraints being imposed to the weight indexed by $\mathcal{J}_2$. Then, we will show
\begin{equation}\label{eq:key1}
  \min_{\bgamma \in \Gamma}R_n(\bgamma | \mathcal{M}, \bmu) = \min_{\bgamma \in \Gamma_2}R_n(\bgamma | \mathcal{M}, \bmu).
\end{equation}
It is obvious that $\min_{\bgamma \in \Gamma}R_n(\bgamma | \mathcal{M}, \bmu) \leq \min_{\bgamma \in \Gamma_2}R_n(\bgamma | \mathcal{M}, \bmu)$ since $\Gamma_2 \subseteq \Gamma$. To prove (\ref{eq:key1}), it suffices to show
\begin{equation}\label{eq:key2}
  \bgamma^*|\mathcal{M}\triangleq \argmin_{\bgamma \in \Gamma}R_n(\bgamma | \mathcal{M}, \bmu) \in \Gamma_2.
\end{equation}
The argument (\ref{eq:key2}) can be established by utilizing a proof by contradiction. Suppose $\bgamma^*|\mathcal{M} \notin \Gamma_2$, then there must exist a $l \in \mathcal{I}_2$ such that $\gamma^*_l > \gamma^*_{l+1}$ and $\|\bmu_{l|\mathcal{M}}\|^2/\sigma^2_{l|\mathcal{M}}< \|\bmu_{l+1|\mathcal{M}}\|^2/\sigma^2_{l+1|\mathcal{M}}$. Define a weight vector $\bgamma^{**}|\mathcal{M}=(\bgamma^{**}_1,\ldots,\bgamma^{**}_{M_n})^{\top}$ such that \begin{equation}\label{eq:gamma_ss_1}
\gamma^{**}_m=\left\{\begin{array}{ll}
\gamma^*_m &\quad m \neq l,l+1, \\
\gamma^*_l &\quad m=l, \\
\gamma^*_l & \quad m=l+1,
\end{array}\right.
\end{equation}
when $\gamma^*_l \leq \|\bmu_{l|\mathcal{M}}\|^2/\sigma^2_{l|\mathcal{M}}$, and
\begin{equation}\label{eq:gamma_ss_2}
\gamma^{**}_m=\left\{\begin{array}{ll}
\gamma^*_m &\quad m \neq l,l+1, \\
\gamma^*_{l+1} \vee \left(\frac{\|\bmu_{l|\mathcal{M}}\|^2}{\sigma^2_{l|\mathcal{M}}}\right) &\quad m=l, \\
\gamma^*_{l+1} & \quad m=l+1,
\end{array}\right.
\end{equation}
when $\gamma^*_l > \|\bmu_{l|\mathcal{M}}\|^2/\sigma^2_{l|\mathcal{M}}$. The definitions (\ref{eq:gamma_ss_1})--(\ref{eq:gamma_ss_2}) ensure $\bgamma^{**}|\mathcal{M}\in \Gamma$. Moreover, according to the monotonicity property of the quadratic function $R_n(\bgamma | \mathcal{M}, \bmu)$, we see
$$
R_n(\bgamma^{**} | \mathcal{M}, \bmu) < R_n(\bgamma^{*} | \mathcal{M}, \bmu),
$$
which contradicts the definition of $\bgamma^{*} | \mathcal{M}$. Thus, the relation (\ref{eq:key1}) is established.

Define a reduced candidate model set
\begin{equation*}\label{eq:addd4}
\begin{split}
   \mathcal{M}_{2} & \triangleq \mathcal{M}_1\setminus \{\mathcal{M}_1(l):l\in  \mathcal{I}_{2}\} \\
     & = \{ k_l: l \neq \mathcal{I}_{2} \},
\end{split}
\end{equation*}
which excludes the candidate models indexed by $\mathcal{I}_{2}$. According to (\ref{eq:risk}), it is evident that
\begin{equation}\label{eq:key3}
  \min_{\bgamma \in \Gamma_2}R_n(\bgamma | \mathcal{M}, \bmu)=\min_{\bw \in \mathcal{W}_{|\mathcal{M}_{2}|}}R_n(\bw | \mathcal{M}_{2}, \bmu),
\end{equation}
that is, adding the additional equality constraint on $\bgamma$ is equivalent to reducing the candidate models in $\mathcal{M}$. Combining (\ref{eq:key1}) with (\ref{eq:key3}), we see
\begin{equation*}
  R_n(\bw^*|\mathcal{M}, \bmu) = R_n(\bw^*|\mathcal{M}_2, \bmu).
\end{equation*}
If $\|\bmu_{l|\mathcal{M}_2}\|^2/\sigma^2_{l|\mathcal{M}_2}\geq \|\bmu_{l+1|\mathcal{M}_2}\|^2/\sigma^2_{l+1|\mathcal{M}_2}$ for $2 \leq l \leq |\mathcal{M}_2|-1$, we have $T=2$ and $\mathcal{M}_T=\mathcal{M}_2$. Then, (\ref{eq:optimal_convex}) is proved. If the monotonicity assumption is still violated for $\mathcal{M}_{2}$, then we repeat the above process again until obtaining a $\mathcal{M}_{T}$, where all elements in $\{ \|\bmu_{m|\mathcal{M}_T}\|^2/\sigma^2_{m|\mathcal{M}_T}\}_{m=2}^{|\mathcal{M}_T|}$ are monotonically non-increasing.

To complete the proof, it remains to show that the set $\mathcal{M}_T$ generated by Algorithm~1 coincides with the Definition (\ref{eq:M_T_def}). We just need to prove that for any set $\mathcal{S}$ satisfying $|\mathcal{S}| \geq |\mathcal{M}_T|$ and $\frac{ \|\bmu_{l|\mathcal{S}}\|^2}{\sigma^2_{l|\mathcal{S}}}\geq  \frac{ \|\bmu_{l+1|\mathcal{S}}\|^2}{\sigma^2_{l+1|\mathcal{S}}},\, l=2,\ldots, |\mathcal{S}|-1$, we have $S = \mathcal{M}_T$. When $T=1$, we have $\mathcal{M}_T = \mathcal{M} =\{k_1,k_2,\ldots,k_{M_n} \}$, and $\mathcal{S} = \mathcal{M}_T$ naturally holds if $|\mathcal{S}|=M_n$. When $T=2$, define $\mathcal{M}_T =\{k_{i_1},\ldots,k_{i_{|\mathcal{M}_T|}} \}$. Suppose there exists a set $\mathcal{S}=\{k_{j_1},\ldots,k_{j_{|\mathcal{S}|}} \}$ such that $\frac{ \|\bmu_{l|\mathcal{S}}\|^2}{\sigma^2_{l|\mathcal{S}}}\geq  \frac{ \|\bmu_{l+1|\mathcal{S}}\|^2}{\sigma^2_{l+1|\mathcal{S}}},\, l=2,\ldots, |\mathcal{S}|-1$ but $\mathcal{S} \neq \mathcal{M}_T$, then based on the construction of $\mathcal{M}_2$ in Algorithm~1, there must exist SNRs in $\mathcal{S}$ that are not monotonically nonincreasing, leading to a contradiction. Similar arguments can be made for the cases $T\geq 3$.

\section{Proof of Corollary~\ref{theo:improv2}}

Due to the condition $k_{M_n}=p_n$, we have $\|\bmu_{M_n+1|\mathcal{M}}\|^2 = 0$ and
\begin{equation*}
\begin{split}
   R_n(\tilde{\bw}^*|\mathcal{M}, \bmu) & =\sum_{m=1}^{M_n}\frac{\|\bmu_{m|\mathcal{M}}\|^2\sigma^2_{m|\mathcal{M}}}{\|\bmu_{m|\mathcal{M}}\|^2+\sigma^2_{m|\mathcal{M}}} = \sum_{m=1}^{M_n}\frac{\|\bmu_{(m)|\mathcal{M}}\|^2\sigma^2_{(m)|\mathcal{M}}}{\|\bmu_{(m)|\mathcal{M}}\|^2+\sigma^2_{(m)|\mathcal{M}}}.\\
\end{split}
\end{equation*}
The optimal MA risk $R_n(\tilde{\bw}^*|\mathcal{M}, \bmu)$ is lower bounded by
\begin{equation}\label{eq:A2_1}
  \begin{split}
     R_n(\tilde{\bw}^*|\mathcal{M}, \bmu)
       & = \sum_{m=1}^{m_n^*}\frac{\|\bmu_{(m)|\mathcal{M}}\|^2\sigma^2_{(m)|\mathcal{M}}}{\|\bmu_{(m)|\mathcal{M}}\|^2+\sigma^2_{(m)|\mathcal{M}}}+ \sum_{m=m_n^*+1}^{M_n}\frac{\|\bmu_{(m)|\mathcal{M}}\|^2\sigma^2_{(m)|\mathcal{M}}}{\|\bmu_{(m)|\mathcal{M}}\|^2+\sigma^2_{(m)|\mathcal{M}}}\\
       & = \sum_{m=1}^{m_n^*}\frac{\sigma^2_{(m)|\mathcal{M}}}{1+\frac{\sigma^2_{(m)|\mathcal{M}}}{\|\bmu_{(m)|\mathcal{M}}\|^2}}+ \sum_{m=m_n^*+1}^{M_n}\frac{\|\bmu_{(m)|\mathcal{M}}\|^2}{\frac{\|\bmu_{(m)|\mathcal{M}}\|^2}{\sigma^2_{(m)|\mathcal{M}}}+1}\\
       & \geq \frac{1}{2}\sum_{m=1}^{m_n^*}\sigma^2_{(m)|\mathcal{M}} + \frac{1}{2} \sum_{m=m_n^*+1}^{M_n}\|\bmu_{(m)|\mathcal{M}}\|^2,
  \end{split}
\end{equation}
where the last inequality follows from $\|\bmu_{(m)|\mathcal{M}}\|^2/\sigma^2_{(m)|\mathcal{M}} \geq 1$ when $1 \leq m \leq m_n^*$ and $\|\bmu_{(m)|\mathcal{M}}\|^2/\sigma^2_{(m)|\mathcal{M}} < 1$ when $m_n^*+1 \leq m \leq M_n$. In addition, from the second step of (\ref{eq:A2_1}), we also see that
\begin{equation}\label{eq:A2_2}
  R_n(\tilde{\bw}^*|\mathcal{M}, \bmu)  \leq \sum_{m=1}^{m_n^*}\sigma^2_{(m)|\mathcal{M}} +  \sum_{m=m_n^*+1}^{M_n}\|\bmu_{(m)|\mathcal{M}}\|^2.
\end{equation}
Combining (\ref{eq:A2_1}) with (\ref{eq:A2_2}), we have
\begin{equation}\label{eq:A2_step1}
\begin{split}
   R_n(\tilde{\bw}^*|\mathcal{M}, \bmu) & \asymp \sum_{m=1}^{m_n^*}\sigma^2_{(m)|\mathcal{M}} +  \sum_{m=m_n^*+1}^{M_n}\|\bmu_{(m)|\mathcal{M}}\|^2 \asymp \sum_{m=1}^{m_n^*}\sigma^2_{(m)|\mathcal{M}},
\end{split}
\end{equation}
where the second approximation follows from the condition $\sum_{m=m_n^*+1}^{M_n}\|\bmu_{(m)|\mathcal{M}}\|^2=O(\sum_{m=1}^{m_n^*}\sigma^2_{(m)|\mathcal{M}})$.

Now we derive an asymptotic expression for $R_n(\bw^*|\mathcal{M}, \bmu)$. Due to $\|\bmu_{L_n+1|\mathcal{M}_T}\|^2 = \|\bmu_{M_n+1|\mathcal{M}}\|^2 = 0$ and $\sigma^2_{1|\mathcal{M}_T} = \sigma^2_{1|\mathcal{M}}$, we have
\begin{equation}\label{eq:A2_3}
\begin{split}
    R_n(\bw^*|\mathcal{M}, \bmu) =\sigma^2_{1|\mathcal{M}}+\sum_{l=2}^{L_n}\frac{\|\bmu_{l|\mathcal{M}_T}\|^2\sigma^2_{l|\mathcal{M}_T}}{\|\bmu_{l|\mathcal{M}_T}\|^2+\sigma^2_{l|\mathcal{M}_T}}.
\end{split}
\end{equation}
The second term on the right side of (\ref{eq:A2_3}) is lower bounded by
\begin{equation}\label{eq:A2_4}
  \begin{split}
     \sum_{l=2}^{L_n}\frac{\|\bmu_{l|\mathcal{M}_T}\|^2\sigma^2_{l|\mathcal{M}_T}}{\|\bmu_{l|\mathcal{M}_T}\|^2+\sigma^2_{l|\mathcal{M}_T}} & = \sum_{l=2}^{l_n^*}\frac{\|\bmu_{l|\mathcal{M}_T}\|^2\sigma^2_{l|\mathcal{M}_T}}{\|\bmu_{l|\mathcal{M}_T}\|^2+\sigma^2_{l|\mathcal{M}_T}} + \sum_{l=l_n^*+1}^{L_n}\frac{\|\bmu_{l|\mathcal{M}_T}\|^2\sigma^2_{l|\mathcal{M}_T}}{\|\bmu_{l|\mathcal{M}_T}\|^2+\sigma^2_{l|\mathcal{M}_T}} \\
       & = \sum_{l=2}^{l_n^*}\frac{\sigma^2_{l|\mathcal{M}_T}}{1+\frac{\sigma^2_{l|\mathcal{M}_T}}{\|\bmu_{l|\mathcal{M}_T}\|^2}} + \sum_{l=l_n^*+1}^{L_n}\frac{\|\bmu_{l|\mathcal{M}_T}\|^2}{\frac{\|\bmu_{l|\mathcal{M}_T}\|^2}{\sigma^2_{l|\mathcal{M}_T}}+1}\\
       & \geq \frac{1}{2}\sum_{l=2}^{l_n^*}\sigma^2_{l|\mathcal{M}_T} + \frac{1}{2}\sum_{l=l_n^*+1}^{L_n}\|\bmu_{l|\mathcal{M}_T}\|^2,
  \end{split}
\end{equation}
where the inequality follows from $\|\bmu_{l|\mathcal{M}_T}\|^2/\sigma^2_{l|\mathcal{M}_T} \geq 1$ when $2 \leq l \leq l_n^*$ and $\|\bmu_{l|\mathcal{M}_T}\|^2/\sigma^2_{l|\mathcal{M}_T} < 1$ when $l_n^*+1 \leq l \leq L_n$. According to the second step in (\ref{eq:A2_4}), we also see that
\begin{equation}\label{eq:A2_5}
  \sum_{l=2}^{L_n}\frac{\|\bmu_{l|\mathcal{M}_T}\|^2\sigma^2_{l|\mathcal{M}_T}}{\|\bmu_{l|\mathcal{M}_T}\|^2+\sigma^2_{l|\mathcal{M}_T}} \leq \sum_{l=2}^{l_n^*}\sigma^2_{l|\mathcal{M}_T} + \sum_{l=l_n^*+1}^{L_n}\|\bmu_{l|\mathcal{M}_T}\|^2.
\end{equation}
Combining (\ref{eq:A2_3})--(\ref{eq:A2_5}), we see $R_n(\bw^*|\mathcal{M}, \bmu)$ has the rate
\begin{equation}\label{eq:A2_step2}
  R_n(\bw^*|\mathcal{M}, \bmu) \asymp \sum_{l=1}^{l_n^*}\sigma^2_{l|\mathcal{M}_T} + \sum_{l=l_n^*+1}^{L_n}\|\bmu_{l|\mathcal{M}_T}\|^2 \gtrsim \sum_{l=1}^{l_n^*}\sigma^2_{l|\mathcal{M}_T}.
\end{equation}
Combining (\ref{eq:A2_step1}) with (\ref{eq:A2_step2}), we have
\begin{equation*}
  \frac{R_n(\tilde{\bw}^*|\mathcal{M}, \bmu)}{R_n(\bw^*|\mathcal{M}, \bmu)} \lesssim \frac{\sum_{m=1}^{m_n^*}\sigma^2_{(m)|\mathcal{M}}}{\sum_{l=1}^{l_n^*}\sigma^2_{l|\mathcal{M}_T}}.
\end{equation*}

\section{Proof of Theorem~\ref{theorem:1}}

The proof of this theorem basically follows the same lines as the proofs in \cite{Cavalier2001penalized, Cavalier2002Sharp} but involves additional complexity since the covariances of $\boldsymbol{\varepsilon}_{m|\mathcal{M}}$, $m=1,\ldots,M_n$ are not identity matrixes.

We lower bound the optimal MA risk by
\begin{equation}\label{eq:first}
\begin{split}
   R_n(\bw^* | \mathcal{M},\bmu)& \geq  R_n(\tilde{\bw}^*|\mathcal{M}, \bmu) \\
   & = \sum_{m=1}^{M_n}\min_{0 \leq \gamma_m \leq 1}\left[ \|\bmu_{m|\mathcal{M}}\|^2 ( 1 - \gamma_m)^2  + \sigma^2_{m|\mathcal{M}}\gamma_m^2 \right]+ \|\bmu_{M_n+1|\mathcal{M}}\|^2\\
   & = \sum_{m=1}^{M_n}\frac{\sigma^2_{m|\mathcal{M}}\|\bmu_{m|\mathcal{M}}\|^2}{\sigma^2_{m|\mathcal{M}}+\|\bmu_{m|\mathcal{M}}\|^2}+ \|\bmu_{M_n+1|\mathcal{M}}\|^2.
\end{split}
\end{equation}
The risk of the Stein-type MA estimator is
\begin{equation}\label{eq:MA_risk}
  \mathbb{E}L_n(\hat{\bw} | \mathcal{M},\bmu)=\sum_{m=1}^{M_n}\mathbb{E}\|\hat{\gamma}_m\by_{m|\mathcal{M}}-\bmu_{m|\mathcal{M}} \|^2+ \|\bmu_{M_n+1|\mathcal{M}}\|^2.
\end{equation}
The task now is to find the upper bounds of $\mathbb{E}\|\hat{\gamma}_m\by_{m|\mathcal{M}}-\bmu_{m|\mathcal{M}} \|^2$, $m=1,\ldots,M_n$, respectively. Following \cite{Cavalier2002Sharp}, we consider two different cases:
\begin{equation}\label{eq:case1}
  \|\bmu_{m|\mathcal{M}} \|^2<\varphi_m\sigma^2_{m|\mathcal{M}}/2,
\end{equation}
and
\begin{equation}\label{eq:case2}
  \|\bmu_{m|\mathcal{M}} \|^2\geq\varphi_m\sigma^2_{m|\mathcal{M}}/2.
\end{equation}

We first construct the upper bound for $\mathbb{E}\|\hat{\gamma}_m\by_{m|\mathcal{M}}-\bmu_{m|\mathcal{M}} \|^2$ under (\ref{eq:case1}). Note that
\begin{equation}\label{eq:imp1}
\begin{split}
   &\mathbb{E}\|\hat{\gamma}_m\by_{m|\mathcal{M}}-\bmu_{m|\mathcal{M}} \|^2 = \mathbb{E}\| (\hat{\gamma}_m - 1)\by_{m|\mathcal{M}} + \boldsymbol{\varepsilon}_{m|\mathcal{M}} \|^2\\
     & = \mathbb{E}\|\boldsymbol{\varepsilon}_{m|\mathcal{M}} \|^2 + 2 \mathbb{E} \left\langle (\hat{\gamma}_m - 1)\by_{m|\mathcal{M}}, \boldsymbol{\varepsilon}_{m|\mathcal{M}} \right\rangle + \mathbb{E} \| (\hat{\gamma}_m - 1)\by_{m|\mathcal{M}} \|^2.
\end{split}
\end{equation}
The first term $\mathbb{E}\|\boldsymbol{\varepsilon}_{m|\mathcal{M}} \|^2 = \sigma^2_{m|\mathcal{M}}$. For the second term, we obtain
\begin{equation}\label{eq:sec1}
  \mathbb{E} \left\langle (\hat{\gamma}_m - 1)\by_{m|\mathcal{M}}, \boldsymbol{\varepsilon}_{m|\mathcal{M}} \right\rangle = \mathbb{E}\left[\sum_{i=1}^{n}(\hat{\gamma}_m - 1)y_{m|\mathcal{M},i}\varepsilon_{m|\mathcal{M},i}\right],
\end{equation}
where $y_{m|\mathcal{M},i}$, $\mu_{m|\mathcal{M},i}$, and $\varepsilon_{m|\mathcal{M},i}$ denote the $i$-th elements of $\by_{m|\mathcal{M}}$, $\bmu_{m|\mathcal{M}}$, and $\boldsymbol{\varepsilon}_{m|\mathcal{M}}$, respectively.

Define $A_{m|\mathcal{M}}$ the event $\{\|\by_{m|\mathcal{M}} \|^2\geq \sigma^2_{m|\mathcal{M}}(1+\varphi_m)\}$.
Based on Lemma 1 of \cite{Liu1994Siegel}, we have
\begin{equation}\label{eq:sec11}
\begin{split}
     &  \mathbb{E}\left[(\hat{\gamma}_m - 1)y_{m|\mathcal{M},i}\varepsilon_{m|\mathcal{M},i}\right] = \cov \left[y_{m|\mathcal{M},i}, (\hat{\gamma}_m - 1)y_{m|\mathcal{M},i}\right]\\
     & = \sum_{j=1}^n\cov(y_{m|\mathcal{M},i}, y_{m|\mathcal{M},j})\mathbb{E}\left(\frac{\partial \hat{\gamma}_m }{\partial y_{m|\mathcal{M},j}}y_{m|\mathcal{M},i}\right)+ \var(y_{m|\mathcal{M},i})\mathbb{E}(\hat{\gamma}_m - 1),\\
\end{split}
\end{equation}
where
\begin{equation}\label{eq:sec111}
  \frac{\partial \hat{\gamma}_m }{\partial y_{m|\mathcal{M},j}} = \frac{2\sigma^2_{m|\mathcal{M}}(1+\varphi_m)y_{m|\mathcal{M},j}I(A_{m|\mathcal{M}})}{\|\by_{m|\mathcal{M}}\|^4}.
\end{equation}
Substituting (\ref{eq:sec1})--(\ref{eq:sec111}) into (\ref{eq:imp1}) yields
\begin{equation}\label{eq:imp2}
\begin{split}
   &\quad\,\mathbb{E}\|\hat{\gamma}_m\by_{m|\mathcal{M}}-\bmu_{m|\mathcal{M}} \|^2 = \sigma^2_{m|\mathcal{M}}\\
   & +\mathbb{E}\left\{\frac{4\sigma^2_{m|\mathcal{M}}(1+\varphi_m)I(A_{m|\mathcal{M}})}{\|\by_{m|\mathcal{M}}\|^4}\sum_{i = 1}^n\sum_{j = 1}^n \cov(y_{m|\mathcal{M},i}, y_{m|\mathcal{M},j}) y_{m|\mathcal{M},i}y_{m|\mathcal{M},j}\right\}\\
   & - 2\mathbb{E}\left[\frac{\sigma^2_{m|\mathcal{M}}(1+\varphi_m)I(A_{m|\mathcal{M}})}{\|\by_{m|\mathcal{M}}\|^2}+I(\bar{A}_{m|\mathcal{M}})\right]\sum_{i=1}^{n}\var(y_{m|\mathcal{M},i})\\
   & + \mathbb{E}\left[\frac{\sigma^4_{m|\mathcal{M}}(1+\varphi_m)^2I(A_{m|\mathcal{M}})}{\|\by_{m|\mathcal{M}}\|^2}+\|\by_{m|\mathcal{M}}\|^2 I(\bar{A}_{m|\mathcal{M}})\right].
\end{split}
\end{equation}
Since $\sum_{i=1}^{n}\var(y_{m|\mathcal{M},i})=\sigma^2_{m|\mathcal{M}}$ and
\begin{equation}
\begin{split}
     &\sum_{i = 1}^n\sum_{j = 1}^n \cov(y_{m|\mathcal{M},i}, y_{m|\mathcal{M},j}) y_{m|\mathcal{M},i}y_{m|\mathcal{M},j} =\sigma^2 \by_{m|\mathcal{M}}^{\top}\bD_{m|\mathcal{M}}\by_{m|\mathcal{M}} \leq \sigma^2\|\by_{m|\mathcal{M}}\|^2,
\end{split}
\end{equation}
we have
\begin{equation}\label{eq:A6}
  \mathbb{E}\|\hat{\gamma}_m\by_{m|\mathcal{M}}-\bmu_{m|\mathcal{M}} \|^2 \leq \|\bmu_{m|\mathcal{M}} \|^2 + \mathbb{E}\left[W(\|\by_{m|\mathcal{M}}\|^2)I(A_{m|\mathcal{M}})\right],
\end{equation}
where $W(x)$ is the function
\begin{equation}
  W(x)=-x +2\sigma^2_{m|\mathcal{M}}+\frac{\sigma^4_{m|\mathcal{M}}\left[4(1+\varphi_m)/(k_m-k_{m-1})-(1-\varphi_m^2)\right]}{x}.
\end{equation}
Then following the proofs in the Proposition 1 of \cite{Cavalier2002Sharp}, we see that under the condition (\ref{eq:case1}),
\begin{equation}\label{eq:conclud1}
\begin{split}
     & \mathbb{E}\|\hat{\gamma}_m\by_{m|\mathcal{M}}-\bmu_{m|\mathcal{M}} \|^2 \\
     & \leq \frac{1}{1-\varphi_m/2}\frac{\sigma^2_{m|\mathcal{M}}\|\bmu_{m|\mathcal{M}}\|^2}{\sigma^2_{m|\mathcal{M}}+\|\bmu_{m|\mathcal{M}}\|^2}+\frac{8\sigma^2_{m|\mathcal{M}}}{k_m-k_{m-1}}\exp\left[  -\frac{(k_m-k_{m-1})\varphi_m^2}{16(1+2\sqrt{\varphi_m})^2}\right].
\end{split}
\end{equation}

Then, we construct the upper bound under (\ref{eq:case2}). From Theorem 6.2 of \cite{Lehmann1983}, it is evident that
\begin{equation}
  \mathbb{E}\|\hat{\gamma}_m\by_{m|\mathcal{M}}-\bmu_{m|\mathcal{M}} \|^2 \leq \mathbb{E}\|\tilde{\gamma}_m\by_{m|\mathcal{M}}-\bmu_{m|\mathcal{M}} \|^2,
\end{equation}
where
\begin{equation}
  \tilde{\gamma}_m = 1 - \frac{\sigma_{m|\mathcal{M}}^2(1+\varphi_m)}{\|\by_{m|\mathcal{M}}\|^2}
\end{equation}
Similar to (\ref{eq:imp1}), we have
\begin{equation}\label{eq:t2i1}
\begin{split}
   &\mathbb{E}\|\tilde{\gamma}_m\by_{m|\mathcal{M}}-\bmu_{m|\mathcal{M}} \|^2 = \mathbb{E}\| (\tilde{\gamma}_m - 1)\by_{m|\mathcal{M}} + \boldsymbol{\varepsilon}_{m|\mathcal{M}} \|^2\\
     & = \mathbb{E}\|\boldsymbol{\varepsilon}_{m|\mathcal{M}} \|^2 + 2 \mathbb{E} \left\langle (\tilde{\gamma}_m - 1)\by_{m|\mathcal{M}}, \boldsymbol{\varepsilon}_{m|\mathcal{M}} \right\rangle + \mathbb{E} \| (\tilde{\gamma}_m - 1)\by_{m|\mathcal{M}} \|^2.
\end{split}
\end{equation}
And the second term of (\ref{eq:t2i1}) is
\begin{equation}\label{eq:sec1}
  \mathbb{E} \left\langle (\tilde{\gamma}_m - 1)\by_{m|\mathcal{M}}, \boldsymbol{\varepsilon}_{m|\mathcal{M}} \right\rangle = \mathbb{E}\left[\sum_{i=1}^{n}(\tilde{\gamma}_m - 1)y_{m|\mathcal{M},i}\varepsilon_{m|\mathcal{M},i}\right].
\end{equation}
Using Lemma 1 of \cite{Liu1994Siegel} again, we have
\begin{equation}
\begin{split}
   &\mathbb{E}\left[(\tilde{\gamma}_m - 1)y_{m|\mathcal{M},i}\varepsilon_{m|\mathcal{M},i}\right]\\
   &= \sigma_{m|\mathcal{M}}^2(1+\varphi_m)\sum_{j=1}^n \cov (y_{m|\mathcal{M},i}, y_{m|\mathcal{M},j})\mathbb{E}\left(\frac{2y_{m|\mathcal{M},i}y_{m|\mathcal{M},j}}{\|\by_{m|\mathcal{M}}\|^4}\right) \\
     & - \sigma_{m|\mathcal{M}}^2(1+\varphi_m)\var(y_{m|\mathcal{M},i}) \mathbb{E}\left(\frac{1}{\|\by_{m|\mathcal{M}}\|^2}\right).
\end{split}
\end{equation}
Therefore, we have
\begin{equation}\label{eq:t2imp}
\begin{split}
    &   \mathbb{E}\|\tilde{\gamma}_m\by_{m|\mathcal{M}}-\bmu_{m|\mathcal{M}} \|^2\\
     &= \sigma^2_{m|\mathcal{M}} +4\sigma_{m|\mathcal{M}}^2(1+\varphi_m)\sum_{i=1}^{n}\sum_{j=1}^n \cov (y_{m|\mathcal{M},i}, y_{m|\mathcal{M},j})\mathbb{E}\left(\frac{y_{m|\mathcal{M},i}y_{m|\mathcal{M},j}}{\|\by_{m|\mathcal{M}}\|^4}\right) \\
     &\quad - 2\sigma_{m|\mathcal{M}}^2(1+\varphi_m)\sum_{i=1}^{n}\var(y_{m|\mathcal{M},i}) \mathbb{E}\left(\frac{1}{\|\by_{m|\mathcal{M}}\|^2}\right) + \mathbb{E}\left[\frac{\|\by_{m|\mathcal{M}}\|^2\sigma_{m|\mathcal{M}}^4(1+\varphi_m)^2}{\|\by_{m|\mathcal{M}}\|^4}\right]\\
     & \leq \sigma^2_{m|\mathcal{M}} +4\sigma^2\sigma_{m|\mathcal{M}}^2(1+\varphi_m)\mathbb{E}\left(\frac{1}{\|\by_{m|\mathcal{M}}\|^2}\right) - 2\sigma_{m|\mathcal{M}}^4(1+\varphi_m)\mathbb{E}\left(\frac{1}{\|\by_{m|\mathcal{M}}\|^2}\right)\\
      &\quad+ \sigma_{m|\mathcal{M}}^4(1+\varphi_m)^2\mathbb{E}\left(\frac{1}{\|\by_{m|\mathcal{M}}\|^2}\right)\\
     &=\sigma^2_{m|\mathcal{M}} - \left[ 1 - \varphi_m^2 - 4(1+\varphi_m)/(k_m-k_{m-1})  \right] \sigma_{m|\mathcal{M}}^4 \mathbb{E}\frac{1}{\|\by_{m|\mathcal{M}}\|^2}.
\end{split}
\end{equation}
Under Assumption~\ref{ass:2}, the second term of (\ref{eq:t2imp}) is negative. Combining with
\begin{equation}\label{eq:fact_conc}
  \mathbb{E}\left(\frac{1}{\|\by_{m|\mathcal{M}}\|^2}\right) \geq \frac{1}{\mathbb{E}\|\by_{m|\mathcal{M}}\|^2} = \frac{1}{\sigma^2_{m|\mathcal{M}}+\|\bmu_{m|\mathcal{M}}\|^2},
\end{equation}
we have
\begin{equation}\label{eq:conclud2}
\begin{split}
    &   \mathbb{E}\|\tilde{\gamma}_m\by_{m|\mathcal{M}}-\bmu_{m|\mathcal{M}} \|^2\\
     &\leq \frac{\sigma^2_{m|\mathcal{M}}\|\bmu_{m|\mathcal{M}}\|^2}{\sigma^2_{m|\mathcal{M}}+\|\bmu_{m|\mathcal{M}}\|^2} \left[ \frac{\sigma^2_{m|\mathcal{M}}+\|\bmu_{m|\mathcal{M}}\|^2}{\|\bmu_{m|\mathcal{M}}\|^2}-\frac{\sigma^2_{m|\mathcal{M}}}{\|\bmu_{m|\mathcal{M}}\|^2}  + \frac{[\varphi_m^2+8/(k_m-k_{m-1})]\sigma^2_{m|\mathcal{M}}}{\|\bmu_{m|\mathcal{M}}\|^2}  \right] \\
     & \leq \frac{\sigma^2_{m|\mathcal{M}}\|\bmu_{m|\mathcal{M}}\|^2}{\sigma^2_{m|\mathcal{M}}+\|\bmu_{m|\mathcal{M}}\|^2} \left[1+ \frac{2[\varphi_m^2+8/(k_m-k_{m-1})]}{\varphi_m^2} \right].
\end{split}
\end{equation}

In view of Assumption~\ref{ass:2}, we have
\begin{equation}
  \frac{1}{1-\varphi_m / 2} \leq 1+\varphi_m \leq 1+\frac{2\left(\varphi_m^2+8/(k_m-k_{m-1})\right)}{\varphi_m}.
\end{equation}
Combining (\ref{eq:conclud1}) with (\ref{eq:conclud2}), we can conclude that
\begin{equation}\label{eq:conclude}
\begin{split}
   \mathbb{E}\|\hat{\gamma}_m\by_{m|\mathcal{M}}-\bmu_{m|\mathcal{M}} \|^2 & \leq  \left\{(1+2\varphi_m+16/[(k_m-k_{m-1})\varphi_m]\right\}\frac{\sigma^2_{m|\mathcal{M}}\|\bmu_{m|\mathcal{M}}\|^2}{\sigma^2_{m|\mathcal{M}}+\|\bmu_{m|\mathcal{M}}\|^2} \\
     & + \frac{8\sigma^2_{m|\mathcal{M}}}{k_m-k_{m-1}}\exp\left[  -\frac{(k_m-k_{m-1})\varphi_m^2}{16(1+2\sqrt{\varphi_m})^2}\right]
\end{split}
\end{equation}
Under Assumption~\ref{ass:1}, summing up (\ref{eq:conclude}) yields
\begin{equation}\label{eq:conclude3}
  \mathbb{E}L_n(\hat{\bw} | \mathcal{M},\bmu) \leq (1+\bar{\varphi})\sum_{m=1}^{M_n}\frac{\sigma^2_{m|\mathcal{M}}\|\bmu_{m|\mathcal{M}}\|^2}{\sigma^2_{m|\mathcal{M}}+\|\bmu_{m|\mathcal{M}}\|^2}+8c_1\sigma^2+ \|\bmu_{M_n+1|\mathcal{M}}\|^2,
\end{equation}
where $\bar{\varphi} = \max_{1\leq m \leq M_n} 2\varphi_m+16/[(k_m-k_{m-1})\varphi_m$. Substituting (\ref{eq:conclude3}) into (\ref{eq:first}), we have proved the theorem.

\section{Proof of Corollary~\ref{cor:classical_stein}}\label{sec:proof_cor_41}

The risk of the Stein estimator (\ref{eq:est_shr}) is
\begin{equation}\label{eq:shr_risk}
  \mathbb{E}L_n(\hat{\bw}_{\text{\scriptsize STE} }|\mathcal{M},\bmu)=\sum_{m=1}^{M_n}\mathbb{E}\|\hat{\gamma}_{\text{\scriptsize STE},m}\by_{m|\mathcal{M}}-\bmu_{m|\mathcal{M}} \|^2+ \|\bmu_{M_n+1|\mathcal{M}}\|^2.
\end{equation}
The main task is to upper bound $\mathbb{E}\|\hat{\gamma}_{\text{\scriptsize STE},m}\by_{m|\mathcal{M}}-\bmu_{m|\mathcal{M}} \|^2$ for $m=1,\ldots,M_n$.

Based on Theorem 6.2 of \cite{Lehmann1983}, we see
\begin{equation*}
  \mathbb{E}\|\hat{\gamma}_{\text{\scriptsize STE},m}\by_{m|\mathcal{M}}-\bmu_{m|\mathcal{M}} \|^2 \leq \mathbb{E}\|\tilde{\gamma}_{\text{\scriptsize STE},m}\by_{m|\mathcal{M}}-\bmu_{m|\mathcal{M}} \|^2,
\end{equation*}
where
\begin{equation}
  \tilde{\gamma}_{\text{\scriptsize STE},m} = 1 - \frac{\sigma_{m|\mathcal{M}}^2}{\|\by_{m|\mathcal{M}}\|^2}
\end{equation}
Similar to (\ref{eq:t2i1}), we write
\begin{equation}\label{eq:t2i1_2}
\begin{split}
   &\mathbb{E}\|\tilde{\gamma}_{\text{\scriptsize STE},m}\by_{m|\mathcal{M}}-\bmu_{m|\mathcal{M}} \|^2 = \mathbb{E}\| (\tilde{\gamma}_{\text{\scriptsize STE},m} - 1)\by_{m|\mathcal{M}} + \boldsymbol{\varepsilon}_{m|\mathcal{M}} \|^2\\
     & = \mathbb{E}\|\boldsymbol{\varepsilon}_{m|\mathcal{M}} \|^2 + 2 \mathbb{E} \left\langle (\tilde{\gamma}_{\text{\scriptsize STE},m} - 1)\by_{m|\mathcal{M}}, \boldsymbol{\varepsilon}_{m|\mathcal{M}} \right\rangle + \mathbb{E} \| (\tilde{\gamma}_{\text{\scriptsize STE},m} - 1)\by_{m|\mathcal{M}} \|^2.
\end{split}
\end{equation}
The second term of (\ref{eq:t2i1_2}) can be written as
\begin{equation*}
  \mathbb{E} \left\langle (\tilde{\gamma}_{\text{\scriptsize STE},m} - 1)\by_{m|\mathcal{M}}, \boldsymbol{\varepsilon}_{m|\mathcal{M}} \right\rangle = \mathbb{E}\left[\sum_{i=1}^{n}(\tilde{\gamma}_{\text{\scriptsize STE},m} - 1)y_{m|\mathcal{M},i}\varepsilon_{m|\mathcal{M},i}\right].
\end{equation*}
Using Lemma 1 of \cite{Liu1994Siegel}, we have
\begin{equation}
\begin{split}
   &\mathbb{E}\left[(\tilde{\gamma}_{\text{\scriptsize STE},m} - 1)y_{m|\mathcal{M},i}\varepsilon_{m|\mathcal{M},i}\right]\\
   &= \sigma_{m|\mathcal{M}}^2\sum_{j=1}^n \cov (y_{m|\mathcal{M},i}, y_{m|\mathcal{M},j})\mathbb{E}\left(\frac{2y_{m|\mathcal{M},i}y_{m|\mathcal{M},j}}{\|\by_{m|\mathcal{M}}\|^4}\right) \\
     & - \sigma_{m|\mathcal{M}}^2\var(y_{m|\mathcal{M},i}) \mathbb{E}\left(\frac{1}{\|\by_{m|\mathcal{M}}\|^2}\right).
\end{split}
\end{equation}
Therefore, we see
\begin{equation}\label{eq:conclud2_add}
\begin{split}
    &   \mathbb{E}\|\tilde{\gamma}_{\text{\scriptsize STE},m}\by_{m|\mathcal{M}}-\bmu_{m|\mathcal{M}} \|^2\\
     &= \sigma^2_{m|\mathcal{M}} +4\sigma_{m|\mathcal{M}}^2\sum_{i=1}^{n}\sum_{j=1}^n \cov (y_{m|\mathcal{M},i}, y_{m|\mathcal{M},j})\mathbb{E}\left(\frac{y_{m|\mathcal{M},i}y_{m|\mathcal{M},j}}{\|\by_{m|\mathcal{M}}\|^4}\right) \\
     &\quad - 2\sigma_{m|\mathcal{M}}^2(1+\varphi_m)\sum_{i=1}^{n}\var(y_{m|\mathcal{M},i}) \mathbb{E}\left(\frac{1}{\|\by_{m|\mathcal{M}}\|^2}\right) + \mathbb{E}\left[\frac{\|\by_{m|\mathcal{M}}\|^2\sigma_{m|\mathcal{M}}^4}{\|\by_{m|\mathcal{M}}\|^4}\right]\\
     & \leq \sigma^2_{m|\mathcal{M}} +4\sigma^2\sigma_{m|\mathcal{M}}^2\mathbb{E}\left(\frac{1}{\|\by_{m|\mathcal{M}}\|^2}\right) - 2\sigma_{m|\mathcal{M}}^4\mathbb{E}\left(\frac{1}{\|\by_{m|\mathcal{M}}\|^2}\right)+ \sigma_{m|\mathcal{M}}^4\mathbb{E}\left(\frac{1}{\|\by_{m|\mathcal{M}}\|^2}\right)\\
     &=\sigma^2_{m|\mathcal{M}} - \left[ 1 -  4/(k_m-k_{m-1})  \right] \sigma_{m|\mathcal{M}}^4 \mathbb{E}\frac{1}{\|\by_{m|\mathcal{M}}\|^2}\\
     & \leq \sigma^2_{m|\mathcal{M}} - \frac{\left[ 1 -  4/(k_m-k_{m-1})  \right] \sigma_{m|\mathcal{M}}^4}{\sigma^2_{m|\mathcal{M}}+\|\bmu_{m|\mathcal{M}}\|^2} \\
     & = \frac{\sigma^2_{m|\mathcal{M}}\|\bmu_{m|\mathcal{M}}\|^2}{\sigma^2_{m|\mathcal{M}}+\|\bmu_{m|\mathcal{M}}\|^2} + \frac{4\sigma^2\sigma^2_{m|\mathcal{M}}}{\sigma^2_{m|\mathcal{M}}+\|\bmu_{m|\mathcal{M}}\|^2}\\
     & \leq \frac{\sigma^2_{m|\mathcal{M}}\|\bmu_{m|\mathcal{M}}\|^2}{\sigma^2_{m|\mathcal{M}}+\|\bmu_{m|\mathcal{M}}\|^2} + 4\sigma^2,
\end{split}
\end{equation}
where the forth step follows from Assumption~\ref{ass:2_add} and the fact (\ref{eq:fact_conc}). Combining (\ref{eq:shr_risk}) with (\ref{eq:conclud2_add}), we obtain
\begin{equation*}
\begin{split}
   \mathbb{E}L_n(\hat{\bw}_{\text{\scriptsize STE} }|\mathcal{M},\bmu) & \leq \sum_{m=1}^{M_n}\frac{\sigma^2_{m|\mathcal{M}}\|\bmu_{m|\mathcal{M}}\|^2}{\sigma^2_{m|\mathcal{M}}+\|\bmu_{m|\mathcal{M}}\|^2} + \|\bmu_{M_n+1|\mathcal{M}}\|^2+4M_n\sigma^2 \\
     & = R_n(\tilde{\bw}^*|\mathcal{M}, \bmu) + 4M_n\sigma^2,
\end{split}
\end{equation*}
where the last equality follows from (\ref{eq:first}).

\section{Proof of Theorem~\ref{theorem:2}}

We first prove that
\begin{equation}\label{eq:ppc}
  R_n(\bw^*|\mathcal{M},\bmu)\leq (1+\zeta_n)R_n(\bw^*|\mathcal{M}_a,\bmu)+k_1\sigma^2.
\end{equation}
Define an $M_n$-dimensional weight vector $\bar{\bw}=(\bar{w}_1,\ldots,\bar{w}_{M_n})^{\top}$, where $\bar{w}_m=\sum_{j=k_{m-1}+1}^{k_m}w_j^*$, $\bar{\gamma}_m=\sum_{j=m}^{M_n}\bar{w}_m$, and $w_j^*$ is the $j$-th element of $\bw^*|\mathcal{M}_a$. According to (\ref{eq:risk}), we have
\begin{equation}\label{eq:ppc2}
\begin{split}
   R_n(\bar{\bw} | \mathcal{M}, \bmu) & = \sum_{m=1}^{M_n}\left[ \|\bmu_{m|\mathcal{M}}\|^2 ( 1 - \bar{\gamma}_m)^2  + \sigma^2_{m|\mathcal{M}}\bar{\gamma}_m^2 \right]\\
     & = \sum_{m=1}^{M_n}\left[ \sum_{j=k_{m-1}+1}^{k_m}\|\bmu_{j|\mathcal{M}_a}\|^2 ( 1 - \bar{\gamma}_m)^2  + \sigma^2_{m|\mathcal{M}}\bar{\gamma}_m^2 \right]\\
     & \leq \sum_{m=1}^{p_n} \|\bmu_{m|\mathcal{M}_a}\|^2 ( 1 - \gamma_m^*)^2  + \sum_{m=1}^{M_n}\sigma^2_{m|\mathcal{M}}\bar{\gamma}_m^2
\end{split}
\end{equation}
where the inequality follows the fact that $\gamma_j^* \leq \bar{\gamma}_m$ when $k_{m-1}+1\leq j \leq k_m$. Note that
\begin{equation}\label{eq:a39}
\begin{split}
   \sum_{m=1}^{M_n}\sigma^2_{m|\mathcal{M}}\bar{\gamma}_m^2 &        \leq k_1\sigma^2+(1+\zeta_n)\sum_{m=2}^{M_n}(k_{m-1}-k_{m-2})\sigma^2\bar{\gamma}_m^2 \\
     & \leq k_1\sigma^2+(1+\zeta_n)\sigma^2\sum_{j=1}^{p_n}(\gamma_j^*)^2,
\end{split}
\end{equation}
where the second inequality is due to $\bar{\gamma}_m\leq \gamma_j^*$ when $k_{m-2}+1\leq j \leq k_{m-1}$.
Substituting (\ref{eq:a39}) into (\ref{eq:ppc2}), we obtain (\ref{eq:ppc}). Then combining the oracle inequality in Theorem~\ref{theorem:1} with (\ref{eq:ppc}), we proof the theorem.

\section{Proof of Corollary~\ref{coro:1}}\label{sec:proof_42}

First, it is obvious that $\mathcal{M}^*$ satisfies Assumption~\ref{ass:3}. We then verify Assumption~\ref{ass:2}. Note that $1/k_1 = 1/\nu_n$. When $2 \leq m \leq M_n$, we have
\begin{equation}
  k_m - k_{m-1} \geq \lfloor  \nu_n \rho_n(1+\rho_n)^{m-1} \rfloor \geq  (1-\rho_n)\nu_n \rho_n(1+\rho_n)^{m-1},
\end{equation}
where the first inequality is due to the definition of $\mathcal{M}^*$, and the second inequality follows from $\lfloor x \rfloor \geq (1 - \rho_n)x$ when $x\geq \rho_n^{-1}$. Thus when $2 \leq m \leq M_n$ and $n$ is large enough, we obtain
\begin{equation}
  \frac{1}{k_m - k_{m-1}} \leq \frac{1}{(1-\rho_n)\nu_n \rho_n} = O\left( \frac{1}{\nu_n \rho_n} \right)
\end{equation}
where the first inequality is due to
\begin{equation}\label{eq:kmm1}
  k_m - k_{m-1} = \lfloor  \nu_n \rho_n(1+\rho_n)^{m-1} \rfloor \geq  (1-\rho_n)\nu_n \rho_n(1+\rho_n)^{m-1}.
\end{equation}
Since $\nu_n\to \infty$ and $\nu_n \rho_n \to \infty$, Assumption~\ref{ass:2} is naturally satisfied when $n$ is large than some fixed integer $N$. Combining with (\ref{eq:bar_phi}), we also see that $\bar{\varphi}\to 0$ as $n \to \infty$.

We now focus on Assumption~\ref{ass:1}. Because $\varphi_m$ is bounded and $\varphi_m=1/(k_m-k_{m-1})^{\tau}$, then there exist a constant $C_1$ such that
\begin{equation}
  \exp\left[  -\frac{(k_m-k_{m-1})\varphi_m^2}{16(1+2\sqrt{\varphi_m})^2}\right] \leq \exp\left[ -C_1 (k_m-k_{m-1})^{1-2\tau} \right].
\end{equation}
When $m=1$, we have $\exp( -C_1 k_1^{1-2\tau} ) \to 0$. When $m\geq 2$, using (\ref{eq:kmm1}), we have
\begin{equation}
  \sum_{m=2}^{M_n}\exp\left[ -C_1 (k_m-k_{m-1})^{1-2\tau} \right] \leq \sum_{m=1}^{\infty}\exp\left[ -C_1C_2 \left(\nu_n \rho_n(1+\rho_n)^{m-1}\right)^{1-2\tau} \right]\to 0,
\end{equation}
which meets Assumption~\ref{ass:1}. Thus, applying Theorem~\ref{theorem:2}, we proof the corollary.

\end{appendix}

\section*{Acknowledgments}

The first version of this paper was completed during the author's visit to the University of Minnesota in 2022. The author would like to thank Professor Yang Li for his unreserved support during the author's visit to the University of Minnesota. The author also thanks Professor Yuhong Yang for helpful comments on an earlier version of this manuscript.

%#####################################################################################################################################
\bigskip
\baselineskip=18pt

\bibliographystyle{apalike}
\bibliography{bibliography}

\end{sloppypar}
\end{document}